\documentstyle{article}

\def\cirk{\,{\raisebox{.3ex}{\tiny $\circ$}}\,}
\def\Ss{{\mbox{$\cal O$}}}
\def\So{{\mbox{${\cal O}_\omega$}}}
\def\Sn{{\mbox{${\cal O}_n$}}}
\def\Sc{{\mbox{$\Delta$}}}
\def\mj{\mbox{\bf 1}}
\def\qed{\hfill $\Box$}
\def\str{\rightarrow}
\def\N{{\mbox{\bf N}}}
\def\Z{{\mbox{\bf Z}}}
\def\pl{\!+\!}
\def\mn{\!-\!}
\def\M{{\mbox{${\cal M}$}}}
\def\Mo{{\mbox{${\cal M}_0$}}}
\def\A{{\mbox{${\cal A}$}}}
\def\B{{\mbox{${\cal B}$}}}
\def\Bn{{\mbox{${\cal B}_n$}}}
\def\fia{\mbox{$\varphi^{{\mbox{\tiny a}}}$}}
\def\gac{\mbox{$\gamma^{{\mbox{\tiny c}}}$}}
\def\Da{{\mbox{${\cal D}_{\cal A}$}}}
\def\Db{{\mbox{${\cal D}_{\cal B}$}}}
\def\Dbn{{\mbox{${\cal D}_{{\cal B}_n}$}}}
\def\Ea{{\mbox{$E_{\cal A}$}}}
\def\Eb{{\mbox{$E_{\cal B}$}}}
\def\Lo{{\mbox{${\cal L}_\omega$}}}

\def\L2n{{\mbox{${\cal L}_{2n}$}}}

\def\Ko{{\mbox{${\cal K}_\omega$}}}
\def\Kn{{\mbox{${\cal K}_n$}}}
\def\K2n{{\mbox{${\cal K}_{2n}$}}}

\def\J{{\mbox{${\cal J}$}}}
\def\Jo{{\mbox{${\cal J}_\omega$}}}

\def\J2n{{\mbox{${\cal J}_{2n}$}}}

\def\dk#1{\mbox{$\lfloor {#1} \rfloor$}}
\def\gk#1{\mbox{$\lceil {#1} \rceil$}}

\begin{document}

\title
{\bf Simplicial Endomorphisms}

\author{
\\[.05cm]
{\sc Kosta Do\v sen}
\\[.05cm]
\\Mathematical Institute, SANU
\\Knez Mihailova 35, P.O. Box 367
\\11001 Belgrade, Serbia
\\email: kosta@mi.sanu.ac.yu
\\[.05cm]}
\date{ }
\maketitle
\begin{abstract}
\noindent The monoids of simplicial endomorphisms, i.e.\ the
monoids of endomorphisms in the simplicial category, are
submonoids of monoids one finds in Temperley-Lieb algebras, and as
the monoids of Temperley-Lieb algebras are linked to situations
where an endofunctor is adjoint to itself, so the monoids of
simplicial endomorphisms are linked to arbitrary adjoint
situations. This link is established through diagrams of the kind
found in Temperley-Lieb algebras. Results about these matters,
which were previously prefigured up to a point, are here surveyed
and reworked. A presentation of monoids of simplicial
endomorphisms by generators and relations has been given a long
time ago. Here a closely related presentation is given, with
completeness proved in a new and self-contained manner.

\vspace{0.3cm}

\noindent{\it Mathematics Subject Classification}
({\it 2000}): 55U10, 20M50, 57M99, 18G30, 18C15, 18A40

\noindent{\it Keywords}: simplicial category, endomorphisms,
presentation by generators and relations,
Temperley-Lieb algebras, adjunction, monads, triples

\end{abstract}

\section{Introduction}

The simplicial category \Sc, whose arrows are all order-preserving
functions from finite ordinals to finite ordinals, plays a central
role in topology. It stands behind many notions of algebraic
topology, and especially those having to do with homology (see
\cite{ML71}, VII.5, and references therein).

For every object $n \geq 0$ of \Sc\ the endomorphisms $f\!:n \str
n$ of \Sc\ make a monoid (i.e.\ semigroup with unit), which has
been considered in semigroup theory under the name \Sn\ (see
\cite{FGJ}). Multiplication in this monoid is composition of
functions, and the unit is the identity function. The goal of the
first part of this paper is to axiomatize, i.e.\ present by
generators and relations, the monoids \Sn\ for every $n$ (in
particular, for every $n \geq 2$, when \Sn\ is not trivial). A
related presentation of \Sn\ has been given a long time ago in
\cite{A62}, and following \cite{A62} a number of related monoids
have been presented and investigated (see \cite{FGJ} and
references therein). Although it is related to this older
presentation, our presentation of \Sn, and our method of proving
its completeness, have some features, mentioned in the next
section, that hopefully make worthwhile their publication. Matters
concerning this presentation will be exposed here in a
self-contained manner, so as to make the whole paper more
self-contained.

The interest of presenting \Sn\ is in the following. The monoids
\Sn\ are submonoids of monoids one finds in Temperley-Lieb
algebras (with vector multiplication), and these algebras play an
important role in knot theory and low-dimensional topology in the
wake of the approach to knot and link invariants through Jones'
polynomial (see \cite{KL}, \cite{L97}, \cite{DP01} and references
therein). The monoids of Temperley-Lieb algebras are monoids of
endomorphisms in categories where an endofunctor is adjoint to
itself (as shown in \cite{DP01}). The monoids \Sn, on the other
hand, arise in adjoint situations in general.

The paper will be organized as follows. In the first part
of the paper (Sections 2-3), the
monoids \Sn\ are presented, and it is proved that they are the
monoids of endomorphisms of \Sc. This proof is based on the presentation
of elements of \Sn\ in a certain normal form.

In the remaining part of the paper, the results of the first part
are put into context. First (Section 4),
it is shown that \Sc\ is isomorphic
to the free monad (or triple) generated by a single object,
and next (Sections 5-8),
it is shown how \Sc, and hence also \Sn, arise in adjunction.
For all of that we exploit the technique of composition
elimination, analogous to the proof-theoretical technique of
cut elimination (see \cite{KD99}). Diagrams of the kind found
in Temperley-Lieb algebras,
which are bound to any adjoint situation,
are introduced in Section 6 under the name
{\it friezes}, and it is shown how they are linked to \Sc\ and \Sn.
At the end of the paper (Section 9), the
monoids of Temperley-Lieb algebras are presented, and it is shown that
the monoids \Sn\ are submonoids of them. We also
consider there the standard presentation of \Sc\ by generators and relations,
and how it is related to adjunction.

The presentation of \Sn\ by generators and relations, and the
proof of its completeness, in the first part of the paper, though
they differ from that of \cite{A62}, cannot be counted as entirely
new results, but the results concerning \Sn\ in the remaining part
of the paper do not seem to have been explicitly registered
before. Many of the results in this remaining, contextual, part of
the paper are however not quite new, and earlier references I know
about will be mentioned at appropriate places in the text. Some
things have been clearly realized before, while others have been
prefigured, more or less clearly. I didn't manage, however, to
find them all put together. Various authors pay attention to
various things. Here an attempt is made to bridge matters, and
present them systematically, without extraneous material and, to a
great extent, in a self-contained manner. What is maybe most
original in this part is the exhibition of the connection between
friezes and order-preserving functions on the set of natural
numbers in Section 6, which serves to link adjunction to \Sc\ and
\Sn.

\section{The monoids \Sn}

The {\it simplicial category} \Sc\ has as objects all finite
ordinals (including 0) and as arrows all order-preserving
functions (see \cite{ML71}, VII.5); namely, for $\varphi \!: n
\str m$ and $i,j \in \{0,\ldots ,n\mn 1\}$, if $i \leq j$, then
$\varphi (i) \leq \varphi (j)$. (The empty function $\O \!: 0 \str
m$ is also order-preserving.) The category \Sc\ is a strict
monoidal category (see \cite{ML71}, VII.1) with the bifunctor $+
\!: \Sc \times \Sc \str \Sc$, which is addition on objects, while
for $\varphi \!: n \str n'$ and $\psi \!: m \str m'$ we have

\[
(\varphi + \psi)(i)=
\left\{
\begin{array}{ll}
\varphi (i) & {\mbox{\rm if }} 0\leq i \leq n\mn 1
\\[.1cm]
n' \!+ \psi (i\mn n) & {\mbox{\rm if }} n \leq i \leq n\pl m\mn 1.
\end{array}
\right.
\]

\noindent The unit object of \Sc\ is the ordinal 0 (which is also
an initial object).

For every object $n$ of \Sc, the endomorphisms $f\!: n \str n$ of
\Sc\ make a monoid with composition and the identity function $\mj
\!: n \str n$. We shall now axiomatize this monoid, i.e.\ present
it by generators and relations.

For every $n \geq 0$, let us designate this monoid by \Sn. Although,
in principle, we allow $n$ here to be lesser than 2, the interesting
monoids \Sn\ must have $n \geq 2$. The generators of \Sn\ are the
{\it right-forking terms} $p^i$ and the {\it left-forking terms}
$q^i$, for every $i$ such that $0 \leq i \leq n\mn 2$, provided $n \geq 2$.
If $n<2$, then we don't have any of these generators.

The right-forking term $p^i$ stands for the {\it right-forking
endomorphism} $\sigma (p^i)\!: {n \str n}$ of \Sc, which satisfies
$\sigma (p^i)(j)=j$ for every $j \in \{0,\ldots ,n\mn 1\}$
different from $i\pl 1$, while $\sigma (p^i)(i\pl 1)=i$.
Diagrammatically, we have

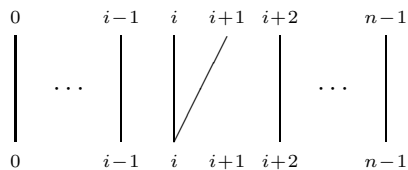
\begin{center}
\begin{picture}(180,60)(0,10)

\put(20,15){\makebox(0,0)[t]{\scriptsize $0$}}
\put(60,15){\makebox(0,0)[t]{\scriptsize $i\mn 1$}}
\put(80,15){\makebox(0,0)[t]{\scriptsize $i$}}
\put(100,15){\makebox(0,0)[t]{\scriptsize $i\pl 1$}}
\put(120,15){\makebox(0,0)[t]{\scriptsize $i\pl 2$}}
\put(160,15){\makebox(0,0)[t]{\scriptsize $n\mn 1$}}

\put(20,65){\makebox(0,0)[b]{\scriptsize $0$}}
\put(60,65){\makebox(0,0)[b]{\scriptsize $i\mn 1$}}
\put(80,65){\makebox(0,0)[b]{\scriptsize $i$}}
\put(100,65){\makebox(0,0)[b]{\scriptsize $i\pl 1$}}
\put(120,65){\makebox(0,0)[b]{\scriptsize $i\pl 2$}}
\put(160,65){\makebox(0,0)[b]{\scriptsize $n\mn 1$}}

\put(40,40){\makebox(0,0){$\cdots$}}
\put(140,40){\makebox(0,0){$\cdots$}}

\put(20,20){\line(0,1){40}}
\put(60,20){\line(0,1){40}}
\put(80,20){\line(0,1){40}}
\put(80,20){\line(1,2){20}}
\put(120,20){\line(0,1){40}}
\put(160,20){\line(0,1){40}}

\end{picture}
\end{center}

\vspace{-.1cm}

The left-forking term $q^i$ stands for the {\it left-forking
endomorphism} $\sigma (q^i)\!: {n \str n}$ of \Sc, which satisfies
$\sigma (q^i)(j)=j$ for every $j \in \{0,\ldots ,n\mn 1\}$
different from $i$, while $\sigma (q^i)(i)=i\pl 1$.
Diagrammatically, we have

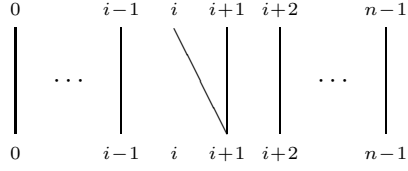
\begin{center}
\begin{picture}(180,60)(0,10)

\put(20,15){\makebox(0,0)[t]{\scriptsize $0$}}
\put(60,15){\makebox(0,0)[t]{\scriptsize $i\mn 1$}}
\put(80,15){\makebox(0,0)[t]{\scriptsize $i$}}
\put(100,15){\makebox(0,0)[t]{\scriptsize $i\pl 1$}}
\put(120,15){\makebox(0,0)[t]{\scriptsize $i\pl 2$}}
\put(160,15){\makebox(0,0)[t]{\scriptsize $n\mn 1$}}

\put(20,65){\makebox(0,0)[b]{\scriptsize $0$}}
\put(60,65){\makebox(0,0)[b]{\scriptsize $i\mn 1$}}
\put(80,65){\makebox(0,0)[b]{\scriptsize $i$}}
\put(100,65){\makebox(0,0)[b]{\scriptsize $i\pl 1$}}
\put(120,65){\makebox(0,0)[b]{\scriptsize $i\pl 2$}}
\put(160,65){\makebox(0,0)[b]{\scriptsize $n\mn 1$}}

\put(40,40){\makebox(0,0){$\cdots$}}
\put(140,40){\makebox(0,0){$\cdots$}}

\put(20,20){\line(0,1){40}}
\put(60,20){\line(0,1){40}}
\put(100,20){\line(0,1){40}}
\put(100,20){\line(-1,2){20}}
\put(120,20){\line(0,1){40}}
\put(160,20){\line(0,1){40}}

\end{picture}
\end{center}

\vspace{-.1cm}

The generators $p^i$ and $q^i$ correspond to the generators
$v_{n-i}$ and $u_i$, respectively, of \cite{A62} (see also
\cite{FGJ}, \S 2). The more complicated indexing of $v_{n-i}$
stresses the left-right duality, on which \cite{A62} relies for
the proof of completeness of the presentation. When comparing the
presentation below with the closely related presentation of
\cite{A62} one should also bear in mind that multiplication, i.e.\
composition, is written there in reverse order, and that instead
of the finite ordinal $n$, which is equal to $\{0,\ldots,n\mn
1\}$, one has $\{1,\ldots,n\}$. The presentation of \cite{A62} is
more economical, while ours is separative, which means that one
could easily obtain from it the presentations of the two
submonoids generated by right-forking terms alone and left-forking
term alone, should one wish to consider these submonoids. The
proof of completeness for our presentation can easily be adapted
to prove the completeness of the presentations of these
submonoids.

The {\it terms} of \Sn\ are obtained from these generators,
together with the special unit term \mj, with the help of the
binary operation of composition $ \cirk $. The unit term \mj\
stands for the identity function $\mj_n \!: n \str n$ of \Sc. For
terms of \Sn\ we use the letters $x$, $y$, $z$, $t$, $\ldots$ ,
$x_1$, $\ldots$

We assume the following equations for \Sn:

\[
\begin{array}{ll}
\makebox[4em][l]{$(1)$} & \makebox[25em][l]{$x \cirk \mj=x,\quad
\mj \cirk x=x$,}
\\[.05cm]
(2) & x\cirk (y\cirk z)=(x\cirk y)\cirk z,
\end{array}
\]
\[
\begin{array}{ll}
\makebox[4em][l]{$(p1)$} & \makebox[25em][l]
{$p^i \cirk p^j=p^j \cirk p^i, \quad {\mbox{\rm for }} j\pl 1<i$,}\\[.05cm]
(p2) & p^i \cirk p^i=p^i,\\[.05cm]
(p3) & p^i \cirk p^{i+1} \cirk p^i=p^{i+1} \cirk p^i \cirk p^{i+1}=
p^i \cirk p^{i+1},
\end{array}
\]
\[
\begin{array}{ll}
\makebox[4em][l]{$(q1)$} & \makebox[25em][l]
{$q^i \cirk q^j=q^j\cirk q^i, \quad {\mbox{\rm for }} j\pl 1<i$,}\\[.05cm]
(q2) & q^i \cirk q^i=q^i,\\[.05cm]
(q3) & q^i \cirk q^{i+1} \cirk q^i=q^{i+1} \cirk q^i \cirk q^{i+1}=
q^{i+1} \cirk q^i,
\end{array}
\]
\[
\begin{array}{ll}
\makebox[4em][l]{$(pq)$} & \makebox[25em][l]
{$p^i \cirk q^j=q^j \cirk p^i, \quad
{\mbox{\rm for }} j<i {\mbox{ \rm  or }} i\pl 1<j$,}
\end{array}
\]
\[
\begin{array}{ll}
\makebox[4em][l]{$(pq1)$} & \makebox[25em][l]{$p^i \cirk q^i=p^i$,}\\[.05cm]
(pq2) & p^i \cirk q^{i+1}=q^{i+1},
\end{array}
\]
\[
\begin{array}{ll}
\makebox[4em][l]{$(qp1)$} & \makebox[25em][l]{$q^i \cirk p^i=q^i$,}\\[.05cm]
(qp2) & q^{i+1} \cirk p^i =p^i.
\end{array}
\]

For $0\leq j\leq i\leq n\mn 2$, the $p\,$-{\it block} $p^{[i,j]}$
is defined as $p^i \cirk p^{i-1} \cirk \ldots \cirk p^{j+1} \cirk
p^j$; we define analogously the $q$-{\it block} $q^{[i,j]}$ as
$q^i \cirk q^{i-1} \cirk \ldots \cirk q^{j+1} \cirk q^j$. The
$p\,$-block $p^{[i,i]}$, which is defined as $p^i$, will be called
{\it singular}, and analogously for $q$-blocks. A {\it block} is
either a $p\,$-block or a $q$-block.

Although the definitions $p\,$-blocks and $q$-blocks are quite
analogous, the two notions are not symmetrical. This asymmetry
will become quite clear in the next section where we consider the
corresponding endomorphisms after the $p$-Points Lemma and the
$q$-Points Lemma. The proof of completeness for the presentation
of \cite{A62}, which otherwise has a similar inspiration as ours,
uses an analogue of our $q$-blocks and an order-reversed
symmetrical notion involving right-forking terms instead of our
$p\,$-blocks. Our asymmetrical proof of completeness given below
hopefully sheds some new light on the matter. One could imagine a
third possibility for this completeness proof, symmetrical as that
of \cite{A62}, but different, since it would be based on
$p\,$-blocks and an order-reversed symmetrical notion involving
left-forking terms instead of our $q$-blocks.

A term of \Sn\ is in {\it normal form} when it is either \mj\ or it is
of the form

\[
r_1^{[i_1,j_1]} \cirk \ldots \cirk r_m^{[i_m,j_m]}
\]

\noindent where $m \geq 1$, for every $k \in \{1,\ldots,m\}$
we have that $r_k$ is either $p$ or $q$, and, in case $m \geq 2$,
for every $k \in \{1,\ldots,m\mn 1\}$

\begin{quote}
if $r_k$ and $r_{k+1}$ are both $p$ or both $q$, then
$i_k < i_{k+1}$ and $j_k < j_{k+1}$;\\
if $r_k$ is $p$ and $r_{k+1}$ is $q$, then $i_k \pl 1 < j_{k+1}$;\\
if $r_k$ is $q$ and $r_{k+1}$ is $p$, then $i_k < j_{k+1}$.
\end{quote}

\noindent (This normal form is inspired by Jones'
normal form of \cite{J83}, \S 4.1.4,
p. 14; see also \cite{BDP02}, \S 1, and \cite{DP01}, \S 10.
It stems ultimately from the
normal form for symmetric groups suggested by \cite{B11}, Note C.)

For the sake of definiteness, we require that in our normal form
all parentheses are associated to the left (but another arrangement
of parentheses would do as well). In the reduction to normal form below
we will not bother about trivial considerations concerning parentheses.
The associativity equation $(2)$ guarantees that we can move them at will,
as it dispensed us from writing parentheses in $(p3)$ and $(q3)$.

That every term of \Sn\ is equal to a term in normal form will be
demonstrated with the help of an alternative formulation of \Sn,
called the {\it block formulation}, which is obtained as follows.
Instead of the terms $p^i$ and $q^i$, we take the $p\,$-blocks and the
$q$-blocks as generators, we generate terms with these generators, \mj\
and \cirk, and to the monoid equations $(1)$ and $(2)$ we add the
following equations:

\vspace{.3cm}

\noindent {\it pp equations}
\[
\begin{array}{lll}
\makebox[6em][l]{$(\mbox{\rm I}pp)$}
& \makebox[10em][l]
{$\mbox{\rm for } k\pl 1<j,$}
& \makebox[15em][l]
{$p^{[i,j]} \cirk p^{[k,l]} =p^{[k,l]} \cirk p^{[i,j]},$}\\[.1cm]
(\mbox{\rm II}pp) & \mbox{\rm for } k\leq j\leq k\pl 1,
& p^{[i,j]} \cirk p^{[k,l]} =p^{[i,l]},
\end{array}
\]
\vspace{-.2cm}
\[
\begin{array}{lll}
\makebox[6em][l]{$(\mbox{\rm III}pp)$}
& \makebox[10em][l]
{$\mbox{\rm for } j<k,$}
& \makebox[15em][l]{ }\\[.05cm]
(\mbox{\rm III.1}pp) & \mbox{\rm for } l\leq j<k\leq i,
& p^{[i,j]} \cirk p^{[k,l]} =p^{[k-1,l]} \cirk p^{[i,j+1]},\\[.05cm]
(\mbox{\rm III.2}pp) & \mbox{\rm for } l\leq j\leq i<k,
& p^{[i,j]} \cirk p^{[k,l]} =p^{[i,l]} \cirk p^{[k,j+1]},\\[.05cm]
(\mbox{\rm III.3}pp) & \mbox{\rm for } j<l\leq k\leq i,
& p^{[i,j]} \cirk p^{[k,l]} =p^{[k-1,j]} \cirk p^{[i,l]},
\end{array}
\]

\vspace{.2cm}

\noindent {\it qq equations}
\[
\begin{array}{ll}
(\mbox{\rm I}qq) & \mbox{\rm for } k\pl 1<j, \quad
q^{[i,j]} \cirk q^{[k,l]} =q^{[k,l]} \cirk q^{[i,j]},\\[.1cm]
(\mbox{\rm II}qq) & \mbox{\rm for } j\leq k\pl 1 \mbox{\rm{ and }}
(l\leq j \mbox{\rm{ or }} k\leq i), \quad q^{[i,j]} \cirk
q^{[k,l]} =q^{[\max(i,k),\min(j,l)]},
\end{array}
\]

\vspace{.2cm}

\noindent {\it pq equations}
\[
\begin{array}{lll}
 & \mbox{\rm for } j\leq k\pl 1\\[.1cm]
(\mbox{\rm I}pq) & \mbox{\rm for } i<l,
& p^{[k,l]} \cirk q^{[i,j]} =q^{[i,j]} \cirk p^{[k,l]},\\[.1cm]
(\mbox{\rm II}pq) & \mbox{\rm for } l\leq i\\[.075cm]
(\mbox{\rm II.1}pq) & \mbox{\rm for } l\pl 1<j<i\leq k,
& p^{[k,l]} \cirk q^{[i,j]} =p^{[j-2,l]} \cirk q^{[i-1,j]}
\cirk p^{[k,i]},\\[.05cm]
(\mbox{\rm II.1.1}pq) & \mbox{\rm for } l\pl 1<j=i\leq k,
& p^{[k,l]} \cirk q^{[i,j]} =p^{[j-2,l]} \cirk p^{[k,i]},\\[.075cm]
(\mbox{\rm II.2}pq) & \mbox{\rm for } j\leq l\pl 1
\mbox{\rm{ and }} j<i\leq k,
& p^{[k,l]} \cirk q^{[i,j]} =q^{[i-1,j]} \cirk p^{[k,i]},\\[.05cm]
(\mbox{\rm II.2.1}pq) & \mbox{\rm for } j\leq l\pl 1
\mbox{\rm{ and }} j=i\leq k,
& p^{[k,l]} \cirk q^{[i,j]} =p^{[k,i]},\\[.075cm]
(\mbox{\rm II.3}pq) & \mbox{\rm for } l\pl 1<j\leq k\pl 1\leq i,
& p^{[k,l]} \cirk q^{[i,j]} =p^{[j-2,l]}\cirk q^{[i,j]},\\[.075cm]
(\mbox{\rm II.4}pq) & \mbox{\rm for } j\leq l\pl 1\leq k\pl 1\leq i,
& p^{[k,l]} \cirk q^{[i,j]} =q^{[i,j]},
\end{array}
\]

\vspace{.2cm}

\noindent {\it qp equations}
\[
\begin{array}{lll}
 & \mbox{\rm for } l\leq i\\[.1cm]
(\mbox{\rm I}qp) & \mbox{\rm for } k\pl 1 <j,
& q^{[i,j]} \cirk p^{[k,l]} =p^{[k,l]} \cirk q^{[i,j]},\\[.1cm]
(\mbox{\rm II}qp) & \mbox{\rm for } j\leq k\pl 1\\[.075cm]
(\mbox{\rm II.1}qp) & \mbox{\rm for } l<j<i<k,
& q^{[i,j]} \cirk p^{[k,l]} =p^{[j-1,l]} \cirk q^{[i,j+1]}
\cirk p^{[k,i+1]},\\[.05cm]
(\mbox{\rm II.1.1}qp) & \mbox{\rm for } l<j=i<k,
& q^{[i,j]} \cirk p^{[k,l]} =p^{[j-1,l]} \cirk p^{[k,i+1]},\\[.075cm]
(\mbox{\rm II.2}qp) & \mbox{\rm for } j\leq l\leq i <k,
& q^{[i,j]} \cirk p^{[k,l]} =q^{[i,j]} \cirk p^{[k,i+1]},\\[.075cm]
(\mbox{\rm II.3}qp) & \mbox{\rm for } l<j<i
\mbox{\rm{ and }} k\leq i,
& q^{[i,j]} \cirk p^{[k,l]} =p^{[j-1,l]}\cirk q^{[i,j+1]},\\[.05cm]
(\mbox{\rm II.3.1}qp) & \mbox{\rm for } l<j=i
\mbox{\rm{ and }} k\leq i,
& q^{[i,j]} \cirk p^{[k,l]} =p^{[j-1,l]},\\[.075cm]
(\mbox{\rm II.4}qp) & \mbox{\rm for } j\leq l\leq k\leq i,
& q^{[i,j]} \cirk p^{[k,l]} =q^{[i,j]}.
\end{array}
\]

We verify first that with $p^i$ defined as the singular block $p^{[i,i]}$
and $q^i$ defined as the singular block $q^{[i,i]}$ the equations $(p1)$,
$(p2)$, $\ldots$ , $(qp2)$ of the old formulation of \Sn\ are instances
of the new equations of the block formulation.

The equation $(p1)$ is $(\mbox{\rm I}pp)$ for $i=j$ and $k=l$; the
equation $(p2)$ is $(\mbox{\rm II}pp)$ for $i=j=k=l$; and the equations
$(p3)$ are obtained from $(\mbox{\rm III.2}pp)$ with $i=j=k\mn 1=l$,
and from $(\mbox{\rm III.3}pp)$ with $i=j\pl 1=k=l$.

The equation $(q1)$ is $(\mbox{\rm I}qq)$ for $i=j$ and $k=l$; the
equation $(q2)$ is $(\mbox{\rm II}qq)$ for $i=j=k=l$; and the equations
$(q3)$ are obtained from $(\mbox{\rm II}qq)$ with $i=j=k\mn 1=l$,
and from $(\mbox{\rm II}qq)$ with $i=j\pl 1=k=l$.

The equation $(pq)$ is obtained from $(\mbox{\rm I}pq)$ and
$(\mbox{\rm I}qp)$ with $i=j$ and $k=l$.
The equation $(pq1)$ is $(\mbox{\rm II.2.1}pq)$ for $i=j=k=l$, and the
equation $(pq2)$ is $(\mbox{\rm II.4}pq)$ for $i=j=k\pl 1=l\pl 1$.
The equation $(qp1)$ is $(\mbox{\rm II.4}qp)$ for $i=j=k=l$, and,
finally, the
equation $(qp2)$ is $(\mbox{\rm II.3.1}qp)$ for $i=j=k\pl 1=l\pl 1$.

We have to verify too that in the block formulation we can deduce the
definitions of the blocks with the terms $p^i$ and $q^i$ replaced
by singular blocks; namely, we have to verify

\[
p^{[i,j]}=p^{[i,i]}\cirk p^{[i-1,i-1]}\cirk \ldots \cirk p^{[j+1,j+1]}
\cirk p^{[j,j]},
\]

\noindent and analogously for $p$ replaced by $q$.
This follows readily from
$(\mbox{\rm II}pp)$ and $(\mbox{\rm II}qq)$. Note that in all that
we have not used $(\mbox{\rm III.1}pp)$ and most of the equations
of $(\mbox{\rm II}pq)$ and $(\mbox{\rm II}qp)$. (We have put these
superfluous equations in the block formulation to facilitate the
proof of the Normal Form Lemma below.)

To finish showing that the block formulation of \Sn\ is equivalent to
the old formulation, we have to verify that with blocks defined via
the terms $p^i$ and $q^i$ we can deduce all the equations of the
block formulation from the old equations. This is a lengthy, though
pretty straightforward exercise. For that exercise it is useful
to establish by induction on $i\mn j$ that in the old formulation
we have the equations

\[
\begin{array}{l}
q^{[i,j]} \cirk q^{[i,j+1]}=q^{[i-1,j]} \cirk q^{[i,j]}=
q^{[i,j]} \cirk q^{[i,j]}=q^{[i,j]},\\[.05cm]
p^{[i,j]} \cirk q^{[i+1,j+1]}=q^{[i+1,j+1]},\\[.05cm]
q^{[i,j]} \cirk p^{[i,j]}=q^{[i,j]}.
\end{array}
\]

\noindent Our presentation of \Sn\ is such as to facilitate the
proof of the equations of the block formulation. This, and the
wish to separate the two submonoids generated by right-forking
terms alone and left-forking term alone, make us keep redundant
equations. To show that the equations $(p1)$ and $(q1)$, as well
as the equations $(p2)$ and $(q2)$, which can easily be derived
from $(pq1)$ and $(qp1)$, are redundant, and to remove
redundancies in the equations $(p3)$ and $(q3)$, one has the
derivations in \cite{A62} (\S 2).

Then we can prove the following lemma.

\vspace{.3cm}

\noindent {\sc Normal Form Lemma}. {\it Every term of
\Sn\ is equal to a term in normal form.}

\vspace{.2cm}

\noindent {\it Proof.} We will present a reduction procedure that
transforms every term $t$ of \Sn\ into a term $t'$ in normal form
such that $t=t'$ in \Sn. In proof-theoretical jargon,
we establish that this procedure is strongly normalizing; namely,
that every sequence of reductions terminates in a term in
normal form. (This proves more than what is needed for the lemma:
it would be enough for us if for every term {\it some} sequence
of reductions terminates in a term in normal form.) Our
procedure starts by translating a given term of \Sn\ into a term
of the block formulation, which is achieved by replacing the
terms $p^i$ and $q^i$ with singular blocks.

Take a term in the block formulation of \Sn, and with $r$ standing for
either $p$ or $q$, let subterms of this term of the following forms
be called {\it redexes}:

\[
\begin{array}{lll}
(rr) & r^{[i,j]} \cirk r^{[k,l]} & {\mbox{\rm {for }}} k\leq i
{\mbox{\rm{ or }}} l \leq j,\\[.05cm]
(pq) & p^{[k,l]} \cirk q^{[i,j]} & {\mbox{\rm {for }}} j\leq k\pl 1,\\[.05cm]
(qp) & q^{[i,j]} \cirk p^{[k,l]} & {\mbox{\rm {for }}} l\leq i,\\[.05cm]
(\mj) & x \cirk \mj, \mj \cirk x.
\end{array}
\]

\noindent A reduction of the {\it rr} sort consists in
replacing a redex of the
form $(rr)$ by the corresponding term on the right-hand side of one
of the {\it pp} equations or {\it qq} equations. We define analogously
reductions of the sorts {\it pq}, {\it qp} and \mj\ via the
{\it pq}, {\it qp} and $(1)$ equations, respectively.

Let the {\it weight} of a block $r^{[i,j]}$ be $i\mn j\pl 2$, and for
a term $t$ in the block formulation let $m_1 \geq 0$ be the sum
of the weights of all the blocks in $t$. For any subterm of $t$ of the
form $r^{[i,j]}\cirk \ldots \cirk s^{[k,l]}$, where $s$, as $r$, stands
for either $p$ or $q$, the block $s^{[k,l]}$ is said to be
{\it confronted} with $r^{[i,j]}$ when $r^{[i,j]}\cirk s^{[k,l]}$ is
a redex ($r$ and $s$ may, of course, both be $p$ or $q$, or one may be
$p$ and the other $q$). For any block $r^{[i,j]}$ in $t$ let
$\rho (r^{[i,j]})$ be the number of blocks of $t$ confronted with
$r^{[i,j]}$, and let $m_2 \geq 0$ be the sum of all the numbers
$\rho (r^{[i,j]})$ for all the blocks $r^{[i,j]}$ in $t$ plus the number
of occurrences of \mj\ in $t$. The {\it complexity} $\mu (t)$ is the
ordered pair $(m_1,m_2)$; such pairs are well-ordered lexicographically.

Then we check that if $t'$ is obtained from $t$ by a reduction,
then $\mu (t')$ is strictly smaller than $\mu (t)$. With reductions
based on $(\mbox{\rm I}pp)$, $(\mbox{\rm I}qq)$, $(\mbox{\rm I}pq)$,
$(\mbox{\rm I}qp)$ and $(1)$, the number $m_2$ diminishes,
while $m_1$ doesn't change. With reductions based on all the remaining
equations of the block formulation of \Sn, except equation (2), the
number $m_1$ diminishes. So, by induction on the complexity $\mu (t)$,
we obtain that every term of the block formulation is equal to a term
without redexes, and it is easy to see that a term of the block formulation
is without redexes if and only if it is in normal form. \qed

\vspace{.3cm}

We will see towards the end of the next section that this normal
form is unique, but we need not establish this uniqueness for
proving the completeness of our presentation of \Sn.

\section{\Sn\ and \Sc}

In this section we will show that \Sn\ is the monoid of
endomorphisms of \Sc\ in $n$.

Let $\varphi \!: n \str m$ be an arrow of \Sc, i.e.\ an
order-preserving function from $n$ to $m$, and for every $j \in
\{0,\ldots ,m\mn 1\}$ let $\varphi^{-1}(j)= \{i \mid 0\leq i\leq
n\mn 1 {\mbox{\rm{ and }}} \varphi (i) =j\}$. When
$\varphi^{-1}(j)$ is empty, $j$ will be called an {\it empty
point} of $\varphi$, when $\varphi^{-1}(j)$ is a singleton, $j$ is
a {\it single point} of $\varphi$, and when $\varphi^{-1}(j)$ has
more than one member, $j$ is a {\it multiple point} of $\varphi$.

An {\it e-m pair} of $\varphi$ is
a pair of numbers $(i,j)$ such that $i,j \in \{0,\ldots ,m\mn 1\}$
with $i<j$, the number
$i$ is an empty point of $\varphi$, the number $j$ is a multiple
point of $\varphi$, and for every $k$ such that $i<k<j$ the number $k$
is a single point of $\varphi$. The {\it m-e pairs} $(i,j)$
of $\varphi$ are defined in the same way save that $i$ is a multiple and
$j$ an empty point of $\varphi$. The e-m and m-e pairs of $\varphi$
are called the {\it critical pairs} of $\varphi$. The {\it weight}
of a critical pair $(i,j)$ is $j\mn i$.

The {\it complexity} $\nu (\varphi)$ is the pair $(n_1,n_2)$ where
$n_1 \geq 0$ is the number of empty points of $\varphi$, and
$n_2 \geq 0$ is the minimal weight among the weights of the
critical pairs of $\varphi$; if there are no critical pairs
of $\varphi$, then $n_2$ is 0. The pairs $(n_1,n_2)$ are well-ordered
lexicographically. Then we can prove the following lemma.

\vspace{.3cm}

\noindent {\sc Surjectivity Lemma}. {\it Every endomorphism
$\varphi \!: n\str n$ of \Sc\ is equal to an endomorphism from $n$
to $n$ of \Sc\ generated from the endomorphisms $\sigma (p^i)\!:{n
\str n}$ and $\sigma (q^i)\!:{n \str n}$ of \Sc, where $i\in
\{0,\ldots , n\mn 1\}$, together with the identity function $\mj_n
\!: n\str n$, with the help of composition.}

\vspace{.2cm}

\noindent {\it Proof.} We proceed by induction on the complexity
$\nu (\varphi)$. If $\nu (\varphi) = (0,0)$, then $\varphi = \mj_n$.

Suppose now that $\nu (\varphi) = (n_1,n_2) \neq (0,0)$. This means that
$n_1 \neq 0$ and $n_2 \neq 0$, and among all the critical pairs of
$\varphi$ take one with the minimal weight $n_2$. Let that pair be
$(i,j)$.

Suppose $(i,j)$ is an e-m pair, and let $k= \min \varphi^{-1}(i\pl
1)$. Note we have either that $i\pl 1=j$, or that $i\pl 1 < j$
with $i\pl 1$ being a single point of $\varphi$. Then we define
the endomorphism $\varphi'\!:n \str n$ of \Sc\ by taking that
$\varphi'(l)=\varphi(l)$ for every $l \in \{0,\ldots ,n\mn 1\}$
different from $k$, and $\varphi'(k)=i$. It is easy to check that
$\varphi=\sigma(q^i) \cirk \varphi'$. If $i\pl 1=j$, then in
$\nu(\varphi')=(n'_1, n'_2)$ we have that $n'_1 <n_1$, while if
$i\pl 1<j$, then $n'_1=n_1$ but $n'_2<n_2$. In both cases
$\nu(\varphi')<\nu(\varphi)$, and we can apply the induction
hypothesis to $\varphi'$.

Suppose, on the other hand, $(i,j)$ is an m-e pair, and let $k=
\max \varphi^{-1}(j-1)$. Then we define $\varphi'\!:n \str n$ by
taking that $\varphi'(l)=\varphi(l)$ for every $l \in \{0,\ldots
,n\mn 1\}$ different from $k$, and $\varphi'(k)=j$; we check that
$\varphi=\sigma(p^{j-1}) \cirk \varphi'$. Then we conclude that
$\nu(\varphi')<\nu(\varphi)$, and we apply the induction
hypothesis to $\varphi'$.\qed

\vspace{.3cm}

Let $\sigma$ be the map from the terms of \Sn\ to the arrows of \Sc\
defined inductively, starting from $\sigma(p^i)$ and $\sigma(q^i)$, by
putting

\[
\begin{array}{l}
\sigma(\mj)=\mj_n,
\\
\sigma(x \cirk y)= \sigma(x) \cirk \sigma(y).
\end{array}
\]

\noindent Then we can check easily that $\sigma$ defines a homomorphism from
\Sn\ to \Sc.

\vspace{.3cm}

\noindent {\sc Soundness Lemma}. {\it If $x=y$ in \Sn, then
$\sigma(x)=\sigma(y)$ in \Sc.}

\vspace{.3cm}

For $\varphi\!: n\str n$ an endomorphism of \Sc, an empty point
$i$ of $\varphi$ will be called a {\it bottom p-point} of
$\varphi$ when $\varphi(i) < i$. For $i$ a multiple point of
$\varphi$, a member $j$ of $\varphi^{-1}(i)$ such that $i\leq
j<\max\varphi^{-1}(i)$ will be called a {\it top p-point} of
$\varphi$.

If the diagram of $\varphi\!: 8\str 8$ is the following one:

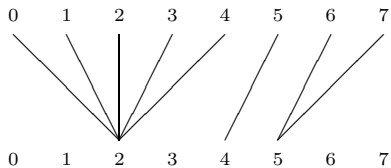
\begin{center}
\begin{picture}(180,60)(0,10)

\put(20,15){\makebox(0,0)[t]{\scriptsize $0$}}
\put(40,15){\makebox(0,0)[t]{\scriptsize $1$}}
\put(60,15){\makebox(0,0)[t]{\scriptsize $2$}}
\put(80,15){\makebox(0,0)[t]{\scriptsize $3$}}
\put(100,15){\makebox(0,0)[t]{\scriptsize $4$}}
\put(120,15){\makebox(0,0)[t]{\scriptsize $5$}}
\put(140,15){\makebox(0,0)[t]{\scriptsize $6$}}
\put(160,15){\makebox(0,0)[t]{\scriptsize $7$}}

\put(20,65){\makebox(0,0)[b]{\scriptsize $0$}}
\put(40,65){\makebox(0,0)[b]{\scriptsize $1$}}
\put(60,65){\makebox(0,0)[b]{\scriptsize $2$}}
\put(80,65){\makebox(0,0)[b]{\scriptsize $3$}}
\put(100,65){\makebox(0,0)[b]{\scriptsize $4$}}
\put(120,65){\makebox(0,0)[b]{\scriptsize $5$}}
\put(140,65){\makebox(0,0)[b]{\scriptsize $6$}}
\put(160,65){\makebox(0,0)[b]{\scriptsize $7$}}

\put(60,20){\line(-1,1){40}}
\put(60,20){\line(-1,2){20}}
\put(60,20){\line(0,1){40}}
\put(60,20){\line(1,2){20}}
\put(60,20){\line(1,1){40}}
\put(100,20){\line(1,2){20}}
\put(120,20){\line(1,2){20}}
\put(120,20){\line(1,1){40}}

\end{picture}
\end{center}

\vspace{-.1cm}

\noindent then 3, 6 and 7 are bottom $p\,$-points and
2, 3 and 6 are top $p\,$-points
of $\varphi$. The denominations ``bottom'' and ``top'' stem from our
putting the domain at the top and the codomain at the bottom of the
diagram. The denomination $p$ is explained by the following lemma.

\vspace{.3cm}

\noindent $p\,$-{\sc{Points Lemma}}. {\it The endomorphism $\sigma(p^{[i,j]})$
has a single bottom $p$-point $i\pl 1$ and a single top $p$-point $j$.}

\vspace{.3cm}

\noindent Diagrammatically, we have for $\sigma(p^{[i,j]})$

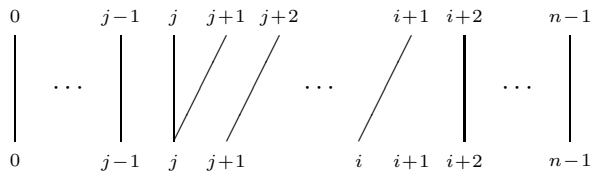
\begin{center}
\begin{picture}(250,60)(0,10)

\put(20,15){\makebox(0,0)[t]{\scriptsize $0$}}
\put(60,15){\makebox(0,0)[t]{\scriptsize $j\mn 1$}}
\put(80,15){\makebox(0,0)[t]{\scriptsize $j$}}
\put(100,15){\makebox(0,0)[t]{\scriptsize $j\pl 1$}}
\put(150,15){\makebox(0,0)[t]{\scriptsize $i$}}
\put(170,15){\makebox(0,0)[t]{\scriptsize $i\pl 1$}}
\put(190,15){\makebox(0,0)[t]{\scriptsize $i\pl 2$}}
\put(230,15){\makebox(0,0)[t]{\scriptsize $n\mn 1$}}

\put(20,65){\makebox(0,0)[b]{\scriptsize $0$}}
\put(60,65){\makebox(0,0)[b]{\scriptsize $j\mn 1$}}
\put(80,65){\makebox(0,0)[b]{\scriptsize $j$}}
\put(100,65){\makebox(0,0)[b]{\scriptsize $j\pl 1$}}
\put(120,65){\makebox(0,0)[b]{\scriptsize $j\pl 2$}}
\put(170,65){\makebox(0,0)[b]{\scriptsize $i\pl 1$}}
\put(190,65){\makebox(0,0)[b]{\scriptsize $i\pl 2$}}
\put(230,65){\makebox(0,0)[b]{\scriptsize $n\mn 1$}}

\put(40,40){\makebox(0,0){$\cdots$}}
\put(135,40){\makebox(0,0){$\cdots$}}
\put(210,40){\makebox(0,0){$\cdots$}}

\put(20,20){\line(0,1){40}}
\put(60,20){\line(0,1){40}}
\put(80,20){\line(0,1){40}}
\put(80,20){\line(1,2){20}}
\put(100,20){\line(1,2){20}}
\put(150,20){\line(1,2){20}}
\put(190,20){\line(0,1){40}}
\put(230,20){\line(0,1){40}}

\end{picture}
\end{center}

\vspace{-.1cm}

\noindent and this proves the lemma.

Let a {\it bottom q-point} of $\varphi$ be a single or multiple
point $i$ of $\varphi$ such that $\min\varphi^{-1}(i)< i$. If $i$
a bottom $q$-point of $\varphi$, then $\min\varphi^{-1}(i)$ will
be called a {\it top q-point} of $\varphi$. These definitions are
explained by the following lemma.

\vspace{.3cm}

\noindent $q$-{\sc{Points Lemma}}. {\it The endomorphism $\sigma(q^{[i,j]})$
has a single bottom $q$-point $i\pl 1$ and a single top $q$-point $j$.}

\vspace{.3cm}

\noindent It is clear from the diagram of $\sigma(q^{[i,j]})$, which looks
as follows:

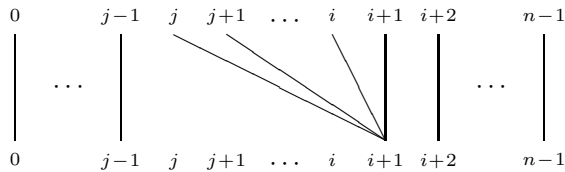
\begin{center}
\begin{picture}(240,60)(0,10)

\put(20,15){\makebox(0,0)[t]{\scriptsize $0$}}
\put(60,15){\makebox(0,0)[t]{\scriptsize $j\mn 1$}}
\put(80,15){\makebox(0,0)[t]{\scriptsize $j$}}
\put(100,15){\makebox(0,0)[t]{\scriptsize $j\pl 1$}}
\put(140,15){\makebox(0,0)[t]{\scriptsize $i$}}
\put(160,15){\makebox(0,0)[t]{\scriptsize $i\pl 1$}}
\put(180,15){\makebox(0,0)[t]{\scriptsize $i\pl 2$}}
\put(220,15){\makebox(0,0)[t]{\scriptsize $n\mn 1$}}

\put(20,65){\makebox(0,0)[b]{\scriptsize $0$}}
\put(60,65){\makebox(0,0)[b]{\scriptsize $j\mn 1$}}
\put(80,65){\makebox(0,0)[b]{\scriptsize $j$}}
\put(100,65){\makebox(0,0)[b]{\scriptsize $j\pl 1$}}
\put(140,65){\makebox(0,0)[b]{\scriptsize $i$}}
\put(160,65){\makebox(0,0)[b]{\scriptsize $i\pl 1$}}
\put(180,65){\makebox(0,0)[b]{\scriptsize $i\pl 2$}}
\put(220,65){\makebox(0,0)[b]{\scriptsize $n\mn 1$}}

\put(40,40){\makebox(0,0){$\cdots$}}
\put(122,11){\makebox(0,0){$\cdots$}}
\put(122,65){\makebox(0,0){$\cdots$}}
\put(200,40){\makebox(0,0){$\cdots$}}

\put(20,20){\line(0,1){40}}
\put(60,20){\line(0,1){40}}
\put(160,20){\line(-2,1){80}}
\put(160,20){\line(-3,2){60}}
\put(160,20){\line(-1,2){20}}
\put(160,20){\line(0,1){40}}
\put(180,20){\line(0,1){40}}
\put(220,20){\line(0,1){40}}

\end{picture}
\end{center}

\vspace{-.1cm}

\noindent that the lemma holds. Note that a bottom $q$-point is not
necessarily a multiple point as in the diagram. (In the next diagram
below, 7 is a bottom $q$-point that is not multiple.)

Consider now the following term of ${\cal S}_{15}$ in normal form

\[
p^{[1,0]}\cirk p^{[2,1]}\cirk p^{[3,3]}\cirk q^{[6,5]}\cirk
q^{[8,6]}\cirk p^{[11,9]}
\]

\noindent which we call $t$. Then the endomorphism $\sigma(t)\!:
15\str 15$ is represented by the diagram

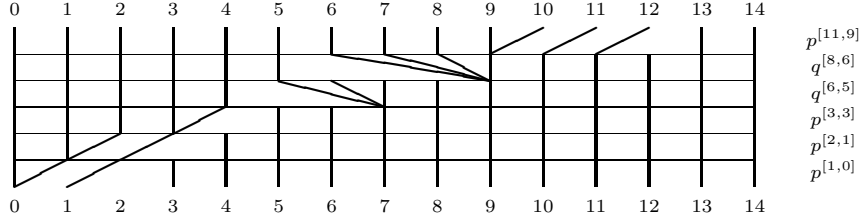
\begin{center}
\begin{picture}(345,80)(0,10)

\thicklines
\put(20,15){\makebox(0,0)[t]{\scriptsize $0$}}
\put(40,15){\makebox(0,0)[t]{\scriptsize $1$}}
\put(60,15){\makebox(0,0)[t]{\scriptsize $2$}}
\put(80,15){\makebox(0,0)[t]{\scriptsize $3$}}
\put(100,15){\makebox(0,0)[t]{\scriptsize $4$}}
\put(120,15){\makebox(0,0)[t]{\scriptsize $5$}}
\put(140,15){\makebox(0,0)[t]{\scriptsize $6$}}
\put(160,15){\makebox(0,0)[t]{\scriptsize $7$}}
\put(180,15){\makebox(0,0)[t]{\scriptsize $8$}}
\put(200,15){\makebox(0,0)[t]{\scriptsize $9$}}
\put(220,15){\makebox(0,0)[t]{\scriptsize $10$}}
\put(240,15){\makebox(0,0)[t]{\scriptsize $11$}}
\put(260,15){\makebox(0,0)[t]{\scriptsize $12$}}
\put(280,15){\makebox(0,0)[t]{\scriptsize $13$}}
\put(300,15){\makebox(0,0)[t]{\scriptsize $14$}}

\put(20,85){\makebox(0,0)[b]{\scriptsize $0$}}
\put(40,85){\makebox(0,0)[b]{\scriptsize $1$}}
\put(60,85){\makebox(0,0)[b]{\scriptsize $2$}}
\put(80,85){\makebox(0,0)[b]{\scriptsize $3$}}
\put(100,85){\makebox(0,0)[b]{\scriptsize $4$}}
\put(120,85){\makebox(0,0)[b]{\scriptsize $5$}}
\put(140,85){\makebox(0,0)[b]{\scriptsize $6$}}
\put(160,85){\makebox(0,0)[b]{\scriptsize $7$}}
\put(180,85){\makebox(0,0)[b]{\scriptsize $8$}}
\put(200,85){\makebox(0,0)[b]{\scriptsize $9$}}
\put(220,85){\makebox(0,0)[b]{\scriptsize $10$}}
\put(240,85){\makebox(0,0)[b]{\scriptsize $11$}}
\put(260,85){\makebox(0,0)[b]{\scriptsize $12$}}
\put(280,85){\makebox(0,0)[b]{\scriptsize $13$}}
\put(300,85){\makebox(0,0)[b]{\scriptsize $14$}}

\put(20,20){\line(0,1){60}}
\put(20,20){\line(2,1){40}}
\put(40,20){\line(2,1){60}}
\put(40,30){\line(0,1){50}}
\put(60,40){\line(0,1){40}}
\put(80,20){\line(0,1){10}}
\put(80,40){\line(0,1){40}}
\put(100,20){\line(0,1){20}}
\put(100,50){\line(0,1){30}}
\put(120,20){\line(0,1){30}}
\put(120,60){\line(0,1){20}}
\put(140,20){\line(0,1){30}}
\put(140,70){\line(0,1){10}}
\put(160,20){\line(0,1){40}}
\put(160,50){\line(-4,1){40}}
\put(160,50){\line(-2,1){20}}
\put(160,70){\line(0,1){10}}
\put(180,20){\line(0,1){40}}
\put(180,70){\line(0,1){10}}
\put(200,20){\line(0,1){60}}
\put(200,60){\line(-6,1){60}}
\put(200,60){\line(-4,1){40}}
\put(200,60){\line(-2,1){20}}
\put(200,70){\line(2,1){20}}
\put(220,20){\line(0,1){50}}
\put(220,70){\line(2,1){20}}
\put(240,20){\line(0,1){50}}
\put(240,70){\line(2,1){20}}
\put(260,20){\line(0,1){50}}
\put(280,20){\line(0,1){60}}
\put(300,20){\line(0,1){60}}

\put(330,30){\makebox(0,0)[t]{\scriptsize $p^{[1,0]}$}}
\put(330,40){\makebox(0,0)[t]{\scriptsize $p^{[2,1]}$}}
\put(330,50){\makebox(0,0)[t]{\scriptsize $p^{[3,3]}$}}
\put(330,60){\makebox(0,0)[t]{\scriptsize $q^{[6,5]}$}}
\put(330,70){\makebox(0,0)[t]{\scriptsize $q^{[8,6]}$}}
\put(330,80){\makebox(0,0)[t]{\scriptsize $p^{[11,9]}$}}

{\linethickness{0.02pt}
\put(20,30){\line(1,0){280}}
\put(20,40){\line(1,0){280}}
\put(20,50){\line(1,0){280}}
\put(20,60){\line(1,0){280}}
\put(20,70){\line(1,0){280}}}

\end{picture}
\end{center}

\vspace{-.1cm}

\noindent which amounts to

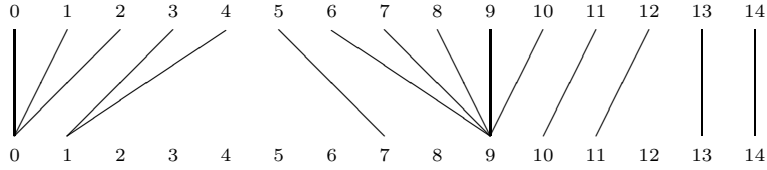
\begin{center}
\begin{picture}(345,60)(0,10)

\put(20,15){\makebox(0,0)[t]{\scriptsize $0$}}
\put(40,15){\makebox(0,0)[t]{\scriptsize $1$}}
\put(60,15){\makebox(0,0)[t]{\scriptsize $2$}}
\put(80,15){\makebox(0,0)[t]{\scriptsize $3$}}
\put(100,15){\makebox(0,0)[t]{\scriptsize $4$}}
\put(120,15){\makebox(0,0)[t]{\scriptsize $5$}}
\put(140,15){\makebox(0,0)[t]{\scriptsize $6$}}
\put(160,15){\makebox(0,0)[t]{\scriptsize $7$}}
\put(180,15){\makebox(0,0)[t]{\scriptsize $8$}}
\put(200,15){\makebox(0,0)[t]{\scriptsize $9$}}
\put(220,15){\makebox(0,0)[t]{\scriptsize $10$}}
\put(240,15){\makebox(0,0)[t]{\scriptsize $11$}}
\put(260,15){\makebox(0,0)[t]{\scriptsize $12$}}
\put(280,15){\makebox(0,0)[t]{\scriptsize $13$}}
\put(300,15){\makebox(0,0)[t]{\scriptsize $14$}}

\put(20,65){\makebox(0,0)[b]{\scriptsize $0$}}
\put(40,65){\makebox(0,0)[b]{\scriptsize $1$}}
\put(60,65){\makebox(0,0)[b]{\scriptsize $2$}}
\put(80,65){\makebox(0,0)[b]{\scriptsize $3$}}
\put(100,65){\makebox(0,0)[b]{\scriptsize $4$}}
\put(120,65){\makebox(0,0)[b]{\scriptsize $5$}}
\put(140,65){\makebox(0,0)[b]{\scriptsize $6$}}
\put(160,65){\makebox(0,0)[b]{\scriptsize $7$}}
\put(180,65){\makebox(0,0)[b]{\scriptsize $8$}}
\put(200,65){\makebox(0,0)[b]{\scriptsize $9$}}
\put(220,65){\makebox(0,0)[b]{\scriptsize $10$}}
\put(240,65){\makebox(0,0)[b]{\scriptsize $11$}}
\put(260,65){\makebox(0,0)[b]{\scriptsize $12$}}
\put(280,65){\makebox(0,0)[b]{\scriptsize $13$}}
\put(300,65){\makebox(0,0)[b]{\scriptsize $14$}}

\put(20,20){\line(0,1){40}}
\put(20,20){\line(1,2){20}}
\put(20,20){\line(1,1){40}}
\put(40,20){\line(1,1){40}}
\put(40,20){\line(3,2){60}}
\put(160,20){\line(-1,1){40}}
\put(200,20){\line(-3,2){60}}
\put(200,20){\line(-1,1){40}}
\put(200,20){\line(-1,2){20}}
\put(200,20){\line(0,1){40}}
\put(200,20){\line(1,2){20}}
\put(220,20){\line(1,2){20}}
\put(240,20){\line(1,2){20}}
\put(280,20){\line(0,1){40}}
\put(300,20){\line(0,1){40}}

\end{picture}
\end{center}

\vspace{-.1cm}

\noindent From this diagram we see that for every $k>11\pl 1$ we have
$\sigma(t)(k)=k$, and for every $k$ such that $9<k \leq 11\pl 1$ we have
$\sigma(t)(k)<k$. If $t'$ is $t$ with $\cirk p^{[12,9]}$ omitted at the
end, then for every $k>8$ we have $\sigma(t')(k)=k$. These findings
are systematized in the following lemma.

\vspace{.3cm}

\noindent {\sc Last Block Lemma}. {\it Let $t$ be the term
$r_1^{[i_1,j_1]} \cirk \ldots \cirk r_m^{[i_m,j_m]}$ of \Sn\ in
normal form. If $r_m$ is $p$ and $i_m\pl 1<k$, or $r_m$ is $q$ and
$i_m<k$, then $\sigma(t)(k)=k$, and if $r_m$ is $p$ and
$j_m<k\leq i_m\pl 1$, then $\sigma(t)(k)<k$.}

\vspace{.2cm}

\noindent {\it Proof.} We proceed by induction on $m$. If $m=1$, then
the lemma is clear from the diagrams of $\sigma(p^{[i,j]})$ and
$\sigma(q^{[i,j]})$ above. Suppose now $m>1$, and the lemma holds for the
term
$r_1^{[i_1,j_1]} \cirk \ldots \cirk r_{m-1}^{[i_{m-1},j_{m-1}]}$,
which we call $t'$.

If $r_m$ is $p$ and $i_m\pl 1<k$, then $\sigma(r_m^{[i_m,j_m]})(k)=k$.
If $r_{m-1}$ is $p$, then $i_{m-1}\pl 1<i_m\pl 1<k$, and
if $r_{m-1}$ is
$q$, then $i_{m-1}<i_{m-1}\pl 1<j_m\pl 1\leq i_m\pl 1<k$.

If $r_m$ is $q$ and $i_m<k$, then $\sigma(r_m^{[i_m,j_m]})(k)=k$.
If $r_{m-1}$ is $p$, then $i_{m-1}\pl 1<j_m\leq i_m<k$,
and if $r_{m-1}$ is
$q$, then $i_{m-1}<i_m<k$.

So we obtain

\[
\begin{array}{ll}
\sigma(t)(k)\!\!\!\!&= \sigma(t')(\sigma(r_m^{[i_m,j_m]})(k))\\
                    &= \sigma(t')(k)\\
                    &= k, \quad {\mbox {\rm {by the induction hypothesis.}}}

\end{array}
\]

If $r_m$ is $p$ and $j_m<k\leq i_m\pl 1$, then
$\sigma(r_m^{[i_m,j_m]})(k)=k\mn 1$.
If $r_{m-1}$ is $p$, then $j_{m-1}<j_m\leq k\mn 1$,
and if $r_{m-1}$ is $q$, then $i_{m-1}<j_m\leq k\mn 1$. Then we
obtain

\[
\begin{array}{ll}
\sigma(t)(k)\!\!\!\!&= \sigma(t')(\sigma(r_m^{[i_m,j_m]})(k))\\
                    &= \sigma(t')(k\mn 1)\\
            &\leq k\mn 1, \quad {\mbox {\rm {by the induction hypothesis}}}\\
                    &< k. \quad \Box
\end{array}
\]

The {\it bottom points} of an endomorphism $\varphi\!: n\str n$ of
\Sc\ are the bottom $p\,$-points and the bottom $q$-points of
$\varphi$, and the {\it top points} of $\varphi$ are the top
$p\,$-points and the top $q$-points of $\varphi$. Note that the
identity function $\mj_n\!: n \str n$ has neither bottom points
nor top points. For terms in normal form other than \mj\ we have
the following lemma.

\vspace{.3cm}

\noindent {\sc Key Lemma}. {\it If $t$ is the term
$r_1^{[i_1,j_1]} \cirk \ldots \cirk r_m^{[i_m,j_m]}$ of \Sn\ in
normal form, then $i_1\pl 1, \ldots , i_m\pl 1$ are all
the bottom points and $j_1, \ldots , j_m$ are all the top points
of $\sigma(t)$. Moreover, $i_k\pl 1$, for $1\leq k\leq m$, is a bottom
$p\,$-point or $q$-point depending on whether $r_k$ is $p$ or $q$, respectively,
and analogously for the top point $j_k$.}

\vspace{.3cm}

Before embarking on the proof of this lemma, we illustrate it with our
diagram

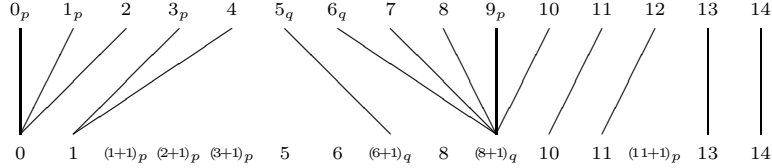
\begin{center}
\begin{picture}(320,60)(0,10)

\put(20,15){\makebox(0,0)[t]{\scriptsize $0$}}
\put(40,15){\makebox(0,0)[t]{\scriptsize $1$}}
\put(60,15){\makebox(0,0)[t]{\tiny $(\!1\!\pl\!1\!)_p$}}
\put(80,15){\makebox(0,0)[t]{\tiny $(\!2\!\pl\!1\!)_p$}}
\put(100,15){\makebox(0,0)[t]{\tiny $(\!3\!\pl\!1\!)_p$}}
\put(120,15){\makebox(0,0)[t]{\scriptsize $5$}}
\put(140,15){\makebox(0,0)[t]{\scriptsize $6$}}
\put(160,15){\makebox(0,0)[t]{\tiny $(\!6\!\pl\!1\!)_q$}}
\put(180,15){\makebox(0,0)[t]{\scriptsize $8$}}
\put(200,15){\makebox(0,0)[t]{\tiny $(\!8\!\pl\!1\!)_q$}}
\put(220,15){\makebox(0,0)[t]{\scriptsize $10$}}
\put(240,15){\makebox(0,0)[t]{\scriptsize $11$}}
\put(260,15){\makebox(0,0)[t]{\tiny $(\!11\!\pl\!1\!)_p$}}
\put(280,15){\makebox(0,0)[t]{\scriptsize $13$}}
\put(300,15){\makebox(0,0)[t]{\scriptsize $14$}}

\put(20,65){\makebox(0,0)[b]{\scriptsize $0_p$}}
\put(40,65){\makebox(0,0)[b]{\scriptsize $1_p$}}
\put(60,65){\makebox(0,0)[b]{\scriptsize $2$}}
\put(80,65){\makebox(0,0)[b]{\scriptsize $3_p$}}
\put(100,65){\makebox(0,0)[b]{\scriptsize $4$}}
\put(120,65){\makebox(0,0)[b]{\scriptsize $5_q$}}
\put(140,65){\makebox(0,0)[b]{\scriptsize $6_q$}}
\put(160,65){\makebox(0,0)[b]{\scriptsize $7$}}
\put(180,65){\makebox(0,0)[b]{\scriptsize $8$}}
\put(200,65){\makebox(0,0)[b]{\scriptsize $9_p$}}
\put(220,65){\makebox(0,0)[b]{\scriptsize $10$}}
\put(240,65){\makebox(0,0)[b]{\scriptsize $11$}}
\put(260,65){\makebox(0,0)[b]{\scriptsize $12$}}
\put(280,65){\makebox(0,0)[b]{\scriptsize $13$}}
\put(300,65){\makebox(0,0)[b]{\scriptsize $14$}}

\put(20,20){\line(0,1){40}}
\put(20,20){\line(1,2){20}}
\put(20,20){\line(1,1){40}}
\put(40,20){\line(1,1){40}}
\put(40,20){\line(3,2){60}}
\put(160,20){\line(-1,1){40}}
\put(200,20){\line(-3,2){60}}
\put(200,20){\line(-1,1){40}}
\put(200,20){\line(-1,2){20}}
\put(200,20){\line(0,1){40}}
\put(200,20){\line(1,2){20}}
\put(220,20){\line(1,2){20}}
\put(240,20){\line(1,2){20}}
\put(280,20){\line(0,1){40}}
\put(300,20){\line(0,1){40}}

\end{picture}
\end{center}

\noindent which corresponds to the term
$p^{[1,0]}\cirk p^{[2,1]}\cirk p^{[3,3]}\cirk q^{[6,5]}\cirk
q^{[8,6]}\cirk p^{[11,9]}$
of ${\cal S}_{15}$ in normal form we had above as an example.
Then $1\pl 1$, $2\pl 1$, $3\pl 1$ and $11\pl 1$ are the bottom
$p\,$-points of this term, $6\pl 1$ and $8\pl 1$ are the bottom
$q$-points, 0, 1, 3 and 9 are the top $p\,$-points, while 5 and 6 are the
top $q$-points.

\vspace{.3cm}

\noindent {\it Proof of the Key Lemma.} We proceed by
induction on $m$. If $m=1$, then we apply the $p\,$-Points Lemma
and the $q$-Points Lemma. Suppose now $m>1$ and the lemma holds for the term
$r_1^{[i_1,j_1]} \cirk \ldots \cirk r_{m-1}^{[i_{m-1},j_{m-1}]}$,
which we call $t'$.

If $r_m$ is $p$, then $i_m\pl 1$ is an empty point of $\sigma(t)$,
because it is an empty point of $\sigma(r_m^{[i_m,j_m]})$ and
a single point of $\sigma(t')$, while $\sigma(t')(i_m\pl 1)= i_m\pl 1$,
by the Last Block Lemma
(if $r_{m-1}$ is $p$, then
$i_{m-1}\pl 1<i_m\pl 1$, and if $r_{m-1}$ is $q$, then
$i_{m-1}<j_m<i_m\pl 1$).
We also have by the Last Block Lemma
that $\sigma(t)(i_m\pl 1)<i_m\pl 1$, and hence $i_m\pl 1$ is a bottom
$p\,$-point of $\sigma(t)$.

If $r_m$ is $p$, then $\sigma(r_m^{[i_m,j_m]})(j_m)=
\sigma(r_m^{[i_m,j_m]})(j_m\pl 1)=j_m$, and hence
$\sigma(t)(j_m)$ $=\sigma(t)(j_m\pl 1)=\sigma(t')(j_m)$.
So both $j_m$ and $j_m\pl 1$
belong to $\sigma(t)^{-1}(\sigma(t')(j_m))$ by using the
Last Block Lemma (if $r_{m-1}$ is $p$, then $j_{m-1}< j_m$, and
if $r_{m-1}$ is $q$, then $i_{m-1}<j_m$). Therefore $j_m$ is a
top $p\,$-point of $\sigma(t)$.

If $r_m$ is $q$, then $j_m=\min\sigma(r_m^{[i_m,j_m]})^{-1}(i_m\pl
1)< i_m\pl 1$. Since $\sigma(t')(i_m\pl 1)=i_m\pl 1$, by the Last
Block Lemma, and $i_m\pl 1$ is a single point of $\sigma(t')$, we
can conclude that $j_m=\min\sigma(t)^{-1}(i_m\pl 1)<i_m\pl 1$.
Hence $i_m\pl 1$ is a bottom $q$-point and $j_m$ is a top
$q$-point of $\sigma(t)$.

It remains to verify that there are no other bottom and top points in
$\sigma(t)$ greater than $i_{m-1}\pl 1$ and $j_{m-1}$, respectively,
save $i_m\pl 1$ and $j_m$. This matter is rather lengthy, but pretty
straightforward. \qed

\vspace{.2cm}

As a consequence of the Key Lemma we have the following two lemmata.

\vspace{.3cm}

\noindent {\sc Auxiliary Lemma}. {\it If $x$ and $y$ are terms
of \Sn\ in normal form and $\sigma(x)=\sigma(y)$ in \Sc, then $x$ and
$y$ are the same term.}

\vspace{.2cm}

\noindent {\it Proof.} If $x$ and $y$ are different terms of \Sn\ in
normal form, then, by the Key Lemma, the functions $\sigma(x)$ and
$\sigma(y)$ must differ with respect to bottom and top points, which
entails that they are different. \qed

\vspace{.3cm}

\noindent {\sc Injectivity Lemma}. {\it If $\sigma(x)=\sigma(y)$ in
\Sc, then $x=y$ in \Sn.}

\vspace{.2cm}

\noindent {\it Proof.} Suppose $\sigma(x)=\sigma(y)$ in \Sc, and let
$x'$ and $y'$ be terms of \Sn\ in normal form such that $x=x'$ and
$y=y'$ in \Sn. Such terms exist according to the Normal Form Lemma
of the preceding section. Then, by the Soundness Lemma,
$\sigma(x)=\sigma(x')$ and $\sigma(y)=\sigma(y')$, and hence
$\sigma(x')=\sigma(y')$. From the Auxiliary Lemma we conclude that
$x'$ is the same term as $y'$, and so $x=y$ in \Sn. \qed

\vspace{.3cm}

The Soundness Lemma and the Injectivity Lemma guarantee that $\sigma$
is a one-one map from \Sn\ to \Sc, and the Surjectivity Lemma guarantees
that $\sigma$ is a map onto all the endomorphisms of \Sc. So \Sn\ is
isomorphic to the monoid of endomorphisms of \Sc\ on the object $n$.

We can also ascertain that for every term $x$ of \Sn\ there is a {\it
unique} term $x'$ in normal form such that $x=x'$ in \Sn.
According to the Normal Form Lemma of the preceding section, take that for
$x'$ and $x''$ in normal form we have $x=x'$ and $x=x''$ in \Sn. Then
we have $x'=x''$ in \Sn, and hence, by the Soundness Lemma,
$\sigma(x')=\sigma(x'')$. Then, by the Auxiliary Lemma,
$x'$ and $x''$ are the same term.

This solves the word problem for \Sn. To check whether $x=y$ in
\Sn\ we could reduce $x$ and $y$ to normal form, according to the
procedure of the proof of the Normal Form Lemma, and then check
whether the normal forms obtained are equal. However, to reduce
a term $x$ of \Sn\ to normal form, now that we have established that
$\sigma$ is an isomorphism, we can proceed more efficiently with
the endomorphism $\sigma(x)$. Just find the bottom and top points
of $\sigma(x)$, from which we immediately obtain $x'$. And to check
whether $x=y$ in \Sn\ it is enough to check whether $\sigma(x)=\sigma(y)$,
which we can do without going via the normal form. But to show that
$\sigma$ is an isomorphism we relied essentially on this normal form.

The fact that $\sigma$ is an isomorphism enables us not only to prove
facts about \Sn\ by going to \Sc, but also facts about \Sc\ by going
to \Sn. For example, we can ascertain that in every endomorphism
$\varphi$ of \Sc\ the number of bottom $p\,$-points is equal to the
number of top $p\,$-points, that the same holds for $q$-points, and
that these points follow each other in a regular manner, as shown by
the normal form. It is not clear how one could prove that directly
in \Sc. We can also ascertain that every endomorphism of \Sc\ is
completely determined by its bottom and top points.

\section{Monads and \Sc}

In this section we will show that \Sc\ is isomorphic to the
free monad (or triple) generated by a single object. This insight
(which perhaps should be traced far back to the {\it Appendice} of
\cite{G58}) may be found
in Lawvere's paper \cite{LAW69} (pp. 148ff; see also \cite{L69}, p. 95,
\cite{A74}, p. 10, and \cite{KD99}, \S 5.9).

A {\it monad} is defined in a standard manner
(see \cite{ML71}, VI) as a quadruple
$\langle \M,T,\eta,\mu\rangle$ where \M\ is a category, $T$ is a
functor from \M\ to \M, while $\eta$ and $\mu$ are natural transformations,
the first from the identity functor on \M\ to $T$
and the second from the composite functor $TT$ to $T$, such that
the following equations hold:

\[
\begin{array}{l}
\mu_a \cirk \eta_{Ta} = \mu_a \cirk T\eta_a = \mj_{Ta},\\
\mu_a \cirk \mu_{Ta} = \mu_a \cirk T\mu_a.
\end{array}
\]

In an alternative, equivalent definition (stemming from \cite{L69}
and akin to a definition of \cite{M76}, \S 1.3, Exercise 12, p.\
32; see \cite{KD99}, \S \S 5.1.5, 5.1.1 and 5.7.3, for the exact
relationship), a monad is a quadruple $\langle \M,T,H,M\rangle$
where

\begin{description}
\item[\hspace{.4cm}] \M\ is a category, \vspace{-.25cm}
\item[\hspace{.4cm}] $T$ is a function from the objects of \M\ to
the objects of \M, \vspace{-.25cm} \item[\hspace{.4cm}] $H$ is a
function that maps the arrows $f\!:a\str b$ of \M\ to the arrows
$HF\!: a\str Tb$ of \M, \vspace{-.25cm} \item[\hspace{.4cm}] $M$
is a function that assigns to every object $b$ of \M\ a function
$M_b$ that maps the arrows $f\!:a\str Tb$ of \M\ to the arrows
$M_bf\!:Ta\str Tb$ of \M,
\end{description}

\noindent and the following equations hold:

\[
\begin{array}{lll}

(H) & Hg\cirk f=H(g\cirk f),
& {\mbox{\rm{  for }}} f\!:a\str b {\mbox{\rm{ and }}} g\!:b\str c,\\
(M) & Mg\cirk Mf=M(Mg\cirk f),
& {\mbox{\rm{  for }}} f\!:a\str Tb {\mbox{\rm{ and }}} g\!:b\str Tc,\\
(H\!M) & Mg\cirk Hf=g\cirk f,
& {\mbox{\rm{  for }}} f\!:a\str b {\mbox{\rm{ and }}} g\!:b\str Tc,\\
(M\!H) & MH\mj_a=\mj_{Ta}

\end{array}
\]

\noindent with indices appropriately assigned to $M$ in $(M)$,
$(H\!M)$ and $(M\!H)$. The two notions of monad are equivalent with the
following definitions:

\[
\begin{array}{ll}
Hf=_{\mbox{\tiny{\rm def}}} \eta_b\cirk f,
& {\mbox{\rm{  for }}} f\!:a\str b,\\
M_b f=_{\mbox{\tiny{\rm def}}} \mu_b\cirk Tf,
& {\mbox{\rm{  for }}} f\!:a\str Tb,\\[.2cm]
Tf=_{\mbox{\tiny{\rm def}}} M_bHf,
& {\mbox{\rm{  for }}} f\!:a\str b,\\
\eta_a=_{\mbox{\tiny{\rm def}}} H\mj_a,\\
\mu_a=_{\mbox{\tiny{\rm def}}} M_a\mj_{Ta}.
\end{array}
\]

The category \Mo\ of the {\it free monad generated by the set of objects}
$\{0\}$ will have as objects finite ordinals. The function $T$ on objects
is defined by $Tn=_{\mbox{\tiny{\rm def}}}n\pl 1$. The
{\it arrow terms} of \Mo, which we call simply {\it terms}, are
defined inductively as follows:

\begin{description}
\item[\hspace{.4cm}] $\mj_n\!:n\str n$ is a term; \vspace{-.25cm}
\item[\hspace{.4cm}] if $f\!:n\str m$ and $g\!:m\str k$ are terms,
then $(g\cirk f)\!:n\str k$ is a term; \vspace{-.25cm}
\item[\hspace{.4cm}] if $f\!:n\str m$ is a term, then $Hf\!:n\str
m\pl 1$ is a term; \vspace{-.25cm} \item[\hspace{.4cm}] if
$f\!:n\str m\pl 1$ is a term, then $Mf\!:n\pl 1\str m\pl 1$ is a
term.
\end{description}

\noindent The expression ``$f\!:n\str m$" is an abbreviation for
``$f$ of type $n\str m$". (The specification of the type $n\str m$
belongs to the metalanguage; in the object language we have only
the terms $f$.) As usual, we omit the outermost parentheses in
$(g\cirk f)$.

Since $T$ is here a one-one function on objects, we don't need to index
$M$ (its index is recovered from the type of $f$ in $Mf$). We impose
the following equations on terms:

\[
\begin{array}{ll}

({\mbox{\it cat}}\:1) & f\cirk \mj_a=f, \quad \mj_b \cirk f=f,
\\[.05cm]
({\mbox{\it cat}}\:2) & h\cirk (g\cirk f)=(h\cirk g)\cirk f,

\end{array}
\]

\noindent and the equations $(H)$, $(M)$, $(H\!M)$ and $(M\!H)$.
(The formal
construction of free monads, in particular those generated by
sets of objects, which may be conceived either as arrowless graphs or as
discrete categories, and the precise sense in which these monads
are free, are explained in \cite{KD99}, \S \S 5.3-6. Formally,
the arrows of \Mo\ are equivalence classes of arrow terms, as the
elements of \Sn\ are equivalence classes of terms of \Sn.)

For the category \Mo\ of the
free monad we have just introduced we can prove the following
proposition due to \cite{L69} (\S 1), which is inspired by Gentzen's
famous cut-elimination technique (see \cite{G35}; see also
\cite{KD99}, \S \S 5.7, 5.8.3).

\vspace{.3cm}

\noindent {\sc Composition Elimination}. {\it For every term $h$ there is
a composition-free term $h'$ such that $h=h'$.}

\vspace{.2cm}

\noindent {\it Proof.} A term of the form $g\cirk f$ where $f$ and $g$
are composition-free is called a {\it topmost composition}. In an arbitrary
term we consider reductions that consist in
replacing a subterm that is a topmost composition
and is of the form on the left-hand side of one of the equations
$({\mbox{\it cat}}\:1)$, $(H)$, $(M)$ and $(H\!M)$ by
the corresponding term on the right-hand side of the equation.

Let the {\it length} of a term be the number of the symbols \mj, \cirk,
$H$ and $M$ in this term (we don't count parentheses). Let the
{\it composition degree} of a term be
the sum of the lengths
of all its subterms of the form $g\cirk f$.
Then it is easy to
check that the length of every topmost composition replaced
in a reduction of the previous paragraph is greater than or equal to the
length of the term by which it is replaced, and that after
every reduction the composition degree of the whole term is strictly
smaller. It remains only to verify that we have covered with
our reductions all possible forms of topmost compositions, and proceed by
induction on the composition degree. \qed

\vspace{.3cm}

Every composition-free term of \Mo\ is of the form $X_n\ldots
X_1\mj_a$, where $n\geq 0$ and $X_i$ is $H$ or $M$. When $a$ is
$0$, then this composition-free term is said to be in {\it normal
form}. Every composition-free term is reduced to a term in normal
form equal to it by replacing $\mj_{n+1}$ with $MH\mj_n$,
according to equation $(M\!H)$. We will see below that this normal
form is unique, i.e.\ that every term of \Mo\ is equal to exactly
one term in normal form.

We will now define inductively a functor $G$ from the category
\Mo\ to the simplicial category \Sc. On objects $G$ is the
identity function, $G(\mj_n)$ is the identity function on $n$, and
$G(g\cirk f)$ is $G(g)\cirk G(f)$. As a set of ordered pairs
$G(Hf)$ is equal to $G(f)$ (but the codomains of these two
functions are different), and for $f\!:n\str m\pl 1$ the set of
ordered pairs $G(Mf)$ is $G(f)\cup \{(n,m)\}$.

We can easily check by induction on the length of derivation that
if $f=g$ in \Mo, then $G(f)=G(g)$ in \Sc. Since it is clear that
$G$ preserves identities and composition, we have that $G$ is a functor
from \Mo\ to \Sc. We will establish below that $G$ is a faithful functor.
Since $G$ is identity on objects, the faithfulness of $G$ amounts to its
being one-one on arrows. We will also establish that $G$ is onto on arrows,
so that we can conclude that the categories \Mo\ and \Sc\ are isomorphic.

\vspace{.3cm}

\noindent {\sc Surjectivity Lemma}. {\it The functor $G$ is onto on arrows.}

\vspace{.2cm}

\noindent {\it Proof.} Take an arbitrary order-preserving function
$\varphi \!: n\str m$. We construct a term $f$ of \Mo\ such that
$G(f)=\varphi$ by induction on $n$.

Suppose $n=0$. Then $\varphi$ is the empty function from $\O$ to $m$.
If $H^m$ stands for a sequence of $m\geq 0$ occurrences of $H$, then
$f$ is $H^m\mj_0$.

Suppose $n=n'\pl 1$, for $n'\geq 0$. In that case $m$ cannot be
$0$; otherwise, $\varphi$ would not exist. Let $\varphi(n')=m'$,
for $m=m'\!+k$, $m'\geq 0$, $k\geq 1$, and let the set of ordered
pairs of the order-preserving function $\varphi'\!:n'\str m'\pl 1$
be defined as the set of ordered pairs of $\varphi$ minus the pair
$(n',m')$. By the induction hypothesis, we have constructed a term
$f'\!:n'\str m'\pl 1$ such that $G(f')=\varphi'$, and $f$ is
$H^{k-1}Mf'$. \qed

\vspace{.3cm}

\noindent {\sc Auxiliary Lemma}. {\it If $f,g\!:k\str l$ are terms
of \Mo\ in normal form and $G(f)=G(g)$, then $f$ is the same term
as $g$.}

\vspace{.2cm}

\noindent {\it Proof.} Suppose $f$, which is $X_n\ldots X_1\mj_0$, and
$g$, which is $Y_m\ldots Y_1\mj_0$, are different terms of \Mo\ of the same
type $k\str l$. Suppose for some $i$ such that $1\leq i\leq n$ and
$1\leq i\leq m$ we have that $X_i$ is different from $Y_i$. Let $j$
be the least such $i$, and let $X_j$ be $H$ while $Y_j$ is $M$.
Let $X_{j-1}\ldots X_1\mj_0$, which is equal to
$Y_{j-1}\ldots Y_1\mj_0$, be of type $r\str s$.

Since $MX_{j-1}\ldots X_1\mj_0$ is defined, we must have that $s=s'\pl 1$
for $s'\geq 0$. Then the pair $(r,s')$ is in
$G(MX_{j-1}\ldots X_1\mj_0)$, and hence also in $G(g)$. But $(r,s')$
cannot belong to $G(f)$, since $HX_{j-1}\ldots X_1\mj_0$ is of type
$r\str s'\pl 2$. So $G(f)\neq G(g)$. If for every $i$ such that
$1\leq i\leq n$
we have that $X_i$ is identical to $Y_i$, then $f$ and $g$ can differ
only if $n<m$. But then $f$ and $g$ cannot be of the same type. \qed

\vspace{.3cm}

\noindent {\sc Injectivity Lemma}. {\it If $f,g\!:k\str l$ are
terms of \Mo\ and $G(f)=G(g)$, then $f=g$ in \Mo.}

\vspace{.2cm}

\noindent {\it Proof.} Suppose $f'$ and $g'$ are the normal forms of
$f$ and $g$ respectively. Since $f=f'$ and $g=g'$ in \Mo, and hence
$G(f)=G(f')$ and $G(g)=G(g')$ in \Sc, from $G(f)=G(g)$ we infer
$G(f')=G(g')$. Then, by the Auxiliary Lemma, $f'$ is the same term
as $g'$, and hence $f=g$ in \Mo. \qed

\vspace{.3cm}

\noindent The proof of this last lemma is analogous to the proof of the
Injectivity Lemma of the preceding section.

As in the preceding section, with the help of
the Auxiliary Lemma and of the functoriality
of $G$ we can ascertain that the normal form of terms of \Mo\ is
unique. For suppose that $f$ and $g$ are normal forms of the same term.
Since both $f$ and $g$ are equal in \Mo\ to this term, we have $f=g$ in
\Mo, and hence $G(f)=G(g)$ in \Sc. Then, by the Auxiliary Lemma,
$f$ and $g$ are the same term.

Since we know that \Sc\ is isomorphic to \Mo, we know that \Sn\
captures the endomorphisms of \Mo. In \Mo\ the endomorphisms that
correspond to the right-forking and left-forking terms of \Sn\ are
defined as follows. If $Tf$ stands for $MHf$, while $T^0$ is the
empty sequence and $T^{n+1}$ is $T^nT$, then we have

\[
\begin{array}{l}
p^i=_{\mbox{\tiny{\rm def}}} T^{n-i-2}HMMH\mj_{T^i0} =
T^{n-i-2}HMT\mj_{T^i0} = T^{n-i-2}(\eta_{TT^i0} \cirk \mu_{T^i0}),
\\[.1cm]
q^i=_{\mbox{\tiny{\rm def}}} T^{n-i-2}MMHH\mj_{T^i0} =
T^{n-i-2}MTH\mj_{T^i0} = T^{n-i-2}(T\eta_{T^i0} \cirk \mu_{T^i0}).
\end{array}
\]

\section{Composition-elimination in adjunction}

An adjunction is defined in a standard manner (see \cite{ML71}, IV)
as a sextuple $\langle \A,\B,F,G,\varphi,\gamma\rangle$
where \A\ and \B\ are categories, $F$ and $G$ are functors, the first
from \B\ to \A\ and the second from \A\ to \B, while $\varphi$ and
$\gamma$ are natural transformations, the first from the composite
functor $FG$ to the identity functor on \A\ and the second from
the identity functor on \B\ to the composite functor $GF$, such that
the following {\it triangular} equations hold:

\[
\quad \varphi_{Fb}\cirk F\gamma_b=\mj_{Fb}, \hspace{2cm}
G\varphi_a\cirk \gamma_{Ga}=\mj_{Ga}.
\]

In an alternative, equivalent definition (see \cite{KD99}, \S \S
4.1-2), an adjunction is a sextuple $\langle \A,\B,F,G,\fia,\gac
\rangle$ where \A, \B, $F$ and $G$ are as in the standard
definition above, while for \fia\ and \gac\ we have that \fia\ is
a function that maps the arrows $f\!:a_1\str a_2$ of \A\ to the
arrows $\fia f\!: FGa_1\str a_2$ of \A, and \gac\ is a function
that maps the arrows $g\!:b_1\str b_2$ of \B\ to the arrows $\gac
g\!: b_1\str GFb_2$ of \B, so that the following equations hold:

\[
\begin{array}{lll}

({\mbox{\it nat}}\:1) \quad & f_2\cirk \fia f_1= \fia (f_2\cirk f_1),
\quad \quad & \gac g_2\cirk g_1=\gac (g_2\cirk g_1),
\\[.05cm]
({\mbox{\it nat}}\:2) \quad & \fia f_2\cirk FGf_1= \fia (f_2\cirk f_1),
\quad & GFg_2\cirk \gac g_1=\gac (g_2\cirk g_1),
\\[.1cm]
(\fia \gac) \quad & \fia f\cirk F\gac g=f\cirk Fg,
\quad & G\fia f\cirk \gac g=Gf\cirk g.

\end{array}
\]

The two notions of adjunction are equivalent with the following definitions:

\[
\begin{array}{ll}

\quad \fia f=_{\mbox{\tiny{\rm def}}} f\cirk \varphi_{a_1},
\hspace{2cm} & \gac g=_{\mbox{\tiny{\rm def}}} \gamma_{b_2}\cirk g,
\\[.05cm]
\quad \varphi_a=_{\mbox{\tiny{\rm def}}} \fia \mj_a,
\hspace{2cm} & \gamma_b=_{\mbox{\tiny{\rm def}}} \gac \mj_b.

\end{array}
\]

We will now describe the {\it free adjunction} $\langle
\A,\B,F,G,\fia,\gac \rangle$ {\it generated by the pair of sets}
$(\O, \{\O\})$, i.e.\ $(0,1)$. The objects of the category \A\ are
generated from $\O$ and those of \B\ from $\{\O\}$. The category
\B\ has as objects words of the form $(GF)^n\O$, where $(GF)^n$
stands for a possibly empty sequence of $n\geq 0$ occurrences of
$GF$. The objects of \A\ are the objects of \B\ with $F$ prefixed.
We use $a$, $a_1$, $\dots$ for the objects of \A\ and $b$, $b_1$,
$\ldots$ for the objects of \B.

The {\it arrow terms} of \A\ and \B, which we call simply {\it terms},
are defined inductively as follows:

\begin{description}

\item[\hspace{.4cm}] $\mj_a\!:a\str a$ is a term of \A;
\vspace{-.25cm} \item[\hspace{.4cm}] $\mj_b\!:b\str b$ is a term
of \B; \vspace{-.05cm} \item[\hspace{.4cm}] if $f_1\!:a_1\str a_2$
and $f_2\!:a_2\str a_3$ are terms of \A, then $f_2\cirk
f_1\!:a_1\str a_3$ is a term of \A; \vspace{-.25cm}
\item[\hspace{.4cm}] if $g_1\!:b_1\str b_2$ and $g_2\!:b_2\str
b_3$ are terms of \B, then $g_2\cirk g_1\!:b_1\str b_3$ is a term
of \B; \vspace{-.05cm} \item[\hspace{.4cm}] if $g\!:b_1\str b_2$
is a term of \B, then $Fg\!:Fb_1\str Fb_2$ is a term of \A;
\vspace{-.25cm} \item[\hspace{.4cm}] if $f\!:a_1\str a_2$ is a
term of \A, then $Gf\!:Ga_1\str Ga_2$ is a term of \B;
\vspace{-.05cm} \item[\hspace{.4cm}] if $f\!:a_1\str a_2$ is a
term of \A, then $\fia f\!:FGa_1\str a_2$ is a term of \A;
\vspace{-.25cm} \item[\hspace{.4cm}] if $g\!:b_1\str b_2$ is a
term of \B, then $\gac g\!:b_1\str GFb_2$ is a term of \B.

\end{description}

We impose on terms the equations $({\mbox{\it cat}}\:1)$,
$({\mbox{\it cat}}\:2)$,

\[
\begin{array}{lll}
({\mbox{\it fun}}\:1) & \quad F\mj_b=\mj_{Fb}, & \quad G\mj_a=\mj_{Ga},
\\[.05cm]
({\mbox{\it fun}}\:2) & \quad Fg_2\cirk Fg_1=F(g_2\cirk g_1), &
\quad Gf_2\cirk Gf_1=G(f_2\cirk f_1),
\end{array}
\]

\noindent and the equations $({\mbox{\it nat}}\:1)$, $({\mbox{\it nat}}\:2)$
and $(\fia \gac)$.
(The formal
construction of free adjunctions, in particular those generated by two
sets of objects, and the precise sense in which these adjunctions
are free, are explained in \cite{KD99}, \S \S 4.2-4. Formally,
the arrows in the categories \A\ and \B\ of our free adjunction are
equivalence classes of arrow terms.)

Our notion of free adjunction is closely related to a
2-categorical notion stemming from \cite{A74} and \cite{SS86}
(which has recently acquired the name ``walking adjunction''; see
\cite{L05}, \S 2.3). Here however we have no need for the
2-categorical context. Our simple notion of free adjunction
suffices to make the connection with the simplicial category \Sc.
This is achieved through a result about the isomorphism of the
categories \B\ and \Sc, towards which we work in the next three
sections of the paper.

For the terms $h$ of the categories \A\ and \B\ of our free
adjunction generated by $(0,1)$ we can prove Composition Elimination
as follows. The statement of this Composition Elimination is as in the
preceding section.

\vspace{.3cm}

\noindent {\it Proof of Composition Elimination.} The notion of
topmost composition is as in the proof of Composition Elimination
in the preceding section. We have reductions for terms which
consist in replacing a subterm that is a topmost composition or is
of the form $F\mj_b$ or $G\mj_a$, according to the equations
$({\mbox{\it cat}}\:1)$, $({\mbox{\it fun}}\:1)$, $({\mbox{\it
fun}}\:2)$, $({\mbox{\it nat}}\:1)$, $({\mbox{\it nat}}\:2)$ and
$(\fia \gac)$, all read from left to right; i.e.\ the subterm
replaced is on the left-hand side, and the term by which it is
replaced is on the right-hand side.

The length of a term is now the number of the symbols \mj, \cirk,
$F$ (applied to arrow terms), $G$ (applied to arrow terms), \fia\
and \gac\ in this term. The composition degree of a term is now defined as
the sum of the lengths
of all its subterms of the form $g\cirk f$ plus the length of the whole
term. Then we reason as in the proof of the preceding section
(see \cite{KD99}, \S \S 4.5, 4.6.3, for details; note that the definition of
composition degree we have here is the one mentioned in parentheses
on p. 118 of \cite{KD99}, and that the other definition mentioned there,
before the parentheses, is not applicable). \qed

\vspace{.3cm}

Every composition-free term of the categories \A\ and \B\ of our
free adjunction is of the
form $X_n\ldots X_1\mj_c$, where $n\geq 0$ and $X_i$ is one of $F$, $G$,
\fia\ and \gac, so that the symbol $F$ can precede immediately only
$G$ or \gac, the symbol $G$ only $F$ or \fia, the symbol \fia\ only
$F$ or \fia, and the symbol \gac\ only $G$ or \gac. When $c$ is $\O$,
this composition-free term is said to be in {\it normal form}.
Every composition-free term is reduced to a term in normal
form equal to it by replacing $\mj_{Fb}$ with $F\mj_b$ and
$\mj_{Ga}$ by $G\mj_a$, according to the
equations $({\mbox{\it fun}}\:1)$ read from right to left.
We will establish in Section 7 that this normal form is unique.

We oriented $({\mbox{\it fun}}\:1)$ from left to right in the
proof of Composition Elimination for the ease of the proof.
Now, for the composition-free normal form, we reverse the
orientation, again to make the matter easier. Our normal form is
an {\it expanded} normal form (as the {\it long} $\beta\eta$
normal forms of the lambda calculus). The
normal form for \Mo\ in the preceding section was also expanded.

\section{Friezes and \Ss}

Let \Ss\ be the monoid of order-preserving endomorphisms of the
set of finite ordinals, i.e.\ the set of natural numbers, \N. In
this monoid multiplication is composition of functions
(order-preserving endomorphisms of \N\ are closed under
composition), and the unit element is the identity function on \N.

In the present section we will consider something we will call
``friezes'', which corresponds to a special kind of
tangle without crossings of knot theory
(see \cite{BZ85}, p. 99, \cite{M96}, Chapter 9, \cite{K95}, Chapter 12).
In \cite{DP01} the term ``frieze'' is used for a different, in
some respects more general notion. We could have called the friezes
of the present paper ``adjunctional friezes'', to distinguish them
from the friezes of \cite{DP01}, but since, in the main body of the paper,
we will have no use for other
friezes save adjunctional ones, we will stick to the shorter
term ``frieze''.
In the next section we will show how the friezes introduced here are
tied to adjunction. In this section we show that our
friezes make a monoid isomorphic to \Ss.

For $M$ an ordered set and for $a,b\in M$ such that
$a<b$, let a {\it segment} $[a,b]$ in $M$
be $\{z\in M\mid a\leq z\leq b\}$. The numbers $a$ and $b$ are the
{\it end points} of $[a,b]$.
We say that $[a,b]$ {\it encloses} $[c,d]$ when $a<c$ and $d<b$.
A set of segments is {\it nonoverlapping} when every two distinct segments
in it are either disjoint or one of these segments encloses the other.
A set $D$ of segments {\it exhausts} $M$ when
all the segments of $D$ are segments in $M$ and
for every $a\in M$ there
is a segment in $D$ one of whose end points is $a$.

A segment $[a,b]$ in $\Z\mn \{0\}$ (i.e.\ the set of integers
without $0$) is called {\it transversal} when $a$ is negative and
$b$ positive; when both of $a$ and $b$ are positive it is a {\it
cup}, and when they are both negative it is a {\it cap}.

A segment in $\Z\mn \{0\}$ is, of course, completely determined
by its end points, and we may as well talk of pairs of integers
$(a,b)$ instead of segments $[a,b]$. We talk of segments to
distinguish them from other sorts of pairs.

A {\it frieze} is a set of nonoverlapping segments
exhausting $\Z\mn \{0\}$ whose
cups are of the form $[2k\pl 2, 2k\pl 3]$ and whose
caps are of the form $[-(2k\pl 2),-(2k\pl 1)]$, for some $k\in \N$.

Friezes may be represented by diagrams. What we do
should be clear from the following example. We draw as follows the
frieze $\{[2,3]$, $[4,5]$, $[10,11]$, $[-2,-1]$,
$[-8,-7]$, $[-3,1]$, $[-4,6]$, $[-5,7]$, $[-6,8]$, $[-9,9]\}\cup
\{[-(k\pl 10),k\pl 12] \mid k\in \N\}$, which, for latter reference,
we call $D_1$:

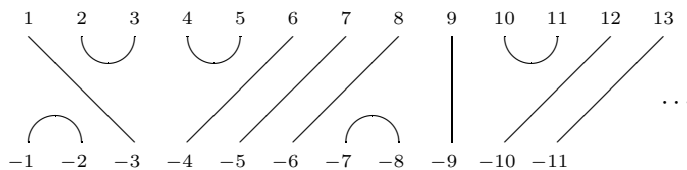
\begin{center}
\begin{picture}(285,60)(0,10)

\put(17,15){\makebox(0,0)[t]{\scriptsize $-1$}}
\put(37,15){\makebox(0,0)[t]{\scriptsize $-2$}}
\put(57,15){\makebox(0,0)[t]{\scriptsize $-3$}}
\put(77,15){\makebox(0,0)[t]{\scriptsize $-4$}}
\put(97,15){\makebox(0,0)[t]{\scriptsize $-5$}}
\put(117,15){\makebox(0,0)[t]{\scriptsize $-6$}}
\put(137,15){\makebox(0,0)[t]{\scriptsize $-7$}}
\put(157,15){\makebox(0,0)[t]{\scriptsize $-8$}}
\put(177,15){\makebox(0,0)[t]{\scriptsize $-9$}}
\put(197,15){\makebox(0,0)[t]{\scriptsize $-10$}}
\put(217,15){\makebox(0,0)[t]{\scriptsize $-11$}}

\put(20,65){\makebox(0,0)[b]{\scriptsize $1$}}
\put(40,65){\makebox(0,0)[b]{\scriptsize $2$}}
\put(60,65){\makebox(0,0)[b]{\scriptsize $3$}}
\put(80,65){\makebox(0,0)[b]{\scriptsize $4$}}
\put(100,65){\makebox(0,0)[b]{\scriptsize $5$}}
\put(120,65){\makebox(0,0)[b]{\scriptsize $6$}}
\put(140,65){\makebox(0,0)[b]{\scriptsize $7$}}
\put(160,65){\makebox(0,0)[b]{\scriptsize $8$}}
\put(180,65){\makebox(0,0)[b]{\scriptsize $9$}}
\put(200,65){\makebox(0,0)[b]{\scriptsize $10$}}
\put(220,65){\makebox(0,0)[b]{\scriptsize $11$}}
\put(240,65){\makebox(0,0)[b]{\scriptsize $12$}}
\put(260,65){\makebox(0,0)[b]{\scriptsize $13$}}

\put(60,20){\line(-1,1){40}}
\put(80,20){\line(1,1){40}}
\put(100,20){\line(1,1){40}}
\put(120,20){\line(1,1){40}}
\put(180,20){\line(0,1){40}}
\put(200,20){\line(1,1){40}}
\put(220,20){\line(1,1){40}}

\put(30,20){\oval(20,20)[t]}
\put(150,20){\oval(20,20)[t]}
\put(50,60){\oval(20,20)[b]}
\put(90,60){\oval(20,20)[b]}
\put(210,60){\oval(20,20)[b]}

\put(265,35){\makebox(0,0){$\cdots$}}

\end{picture}
\end{center}

\vspace{-.1cm}

\noindent This diagram explains the terminology of ``cups'', ``caps'' and
``transversal'' segments.

Note that in a frieze $1$ must be the end point of a transversal
segment. Note also that an adjunctional frieze is uniquely
determined by its transversal segments: the cups and caps need not
be mentioned; but we may as well identify a frieze by its cups and
caps: the transversal segments need not be mentioned.

The {\it unit} frieze \mj\ is $\{[-(k\pl 1), k\pl 1] \mid
k\in \N \}$; in this frieze there are no cups and caps.
It is not so simple to define formally
the composition of two friezes, but it is easy to get
the idea from the following example. If $D_1$ is the
frieze we had above as an example, and $D_2$ is the
frieze that corresponds to the following diagram:

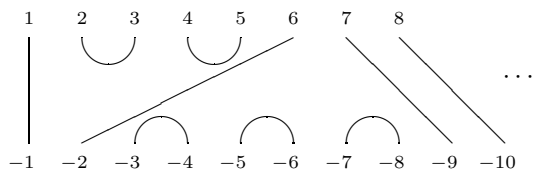
\begin{center}
\begin{picture}(225,60)(0,10)

\put(20,65){\makebox(0,0)[b]{\scriptsize $1$}}
\put(40,65){\makebox(0,0)[b]{\scriptsize $2$}}
\put(60,65){\makebox(0,0)[b]{\scriptsize $3$}}
\put(80,65){\makebox(0,0)[b]{\scriptsize $4$}}
\put(100,65){\makebox(0,0)[b]{\scriptsize $5$}}
\put(120,65){\makebox(0,0)[b]{\scriptsize $6$}}
\put(140,65){\makebox(0,0)[b]{\scriptsize $7$}}
\put(160,65){\makebox(0,0)[b]{\scriptsize $8$}}

\put(17,15){\makebox(0,0)[t]{\scriptsize $-1$}}
\put(37,15){\makebox(0,0)[t]{\scriptsize $-2$}}
\put(57,15){\makebox(0,0)[t]{\scriptsize $-3$}}
\put(77,15){\makebox(0,0)[t]{\scriptsize $-4$}}
\put(97,15){\makebox(0,0)[t]{\scriptsize $-5$}}
\put(117,15){\makebox(0,0)[t]{\scriptsize $-6$}}
\put(137,15){\makebox(0,0)[t]{\scriptsize $-7$}}
\put(157,15){\makebox(0,0)[t]{\scriptsize $-8$}}
\put(177,15){\makebox(0,0)[t]{\scriptsize $-9$}}
\put(197,15){\makebox(0,0)[t]{\scriptsize $-10$}}

\put(20,20){\line(0,1){40}}
\put(40,20){\line(2,1){80}}
\put(180,20){\line(-1,1){40}}
\put(200,20){\line(-1,1){40}}

\put(70,20){\oval(20,20)[t]}
\put(110,20){\oval(20,20)[t]}
\put(150,20){\oval(20,20)[t]}
\put(50,60){\oval(20,20)[b]}
\put(90,60){\oval(20,20)[b]}

\put(205,45){\makebox(0,0){$\cdots$}}

\end{picture}
\end{center}

\vspace{-.1cm}

\noindent then their composition $D_2\cirk D_1$ corresponds to the diagram
obtained by putting the diagram of $D_2$ below
the diagram of $D_1$ in the following manner:

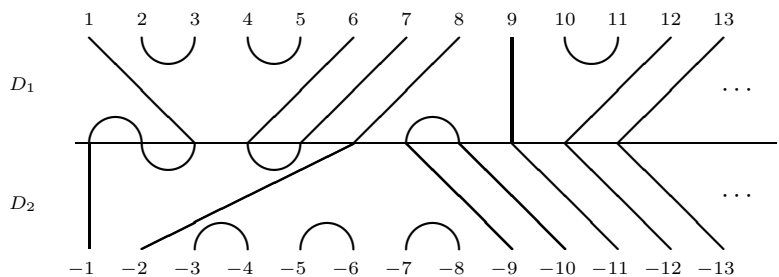
\begin{center}
\begin{picture}(290,100)(0,10)

\put(20,105){\makebox(0,0)[b]{\scriptsize $1$}}
\put(40,105){\makebox(0,0)[b]{\scriptsize $2$}}
\put(60,105){\makebox(0,0)[b]{\scriptsize $3$}}
\put(80,105){\makebox(0,0)[b]{\scriptsize $4$}}
\put(100,105){\makebox(0,0)[b]{\scriptsize $5$}}
\put(120,105){\makebox(0,0)[b]{\scriptsize $6$}}
\put(140,105){\makebox(0,0)[b]{\scriptsize $7$}}
\put(160,105){\makebox(0,0)[b]{\scriptsize $8$}}
\put(180,105){\makebox(0,0)[b]{\scriptsize $9$}}
\put(200,105){\makebox(0,0)[b]{\scriptsize $10$}}
\put(220,105){\makebox(0,0)[b]{\scriptsize $11$}}
\put(240,105){\makebox(0,0)[b]{\scriptsize $12$}}
\put(260,105){\makebox(0,0)[b]{\scriptsize $13$}}

\put(17,15){\makebox(0,0)[t]{\scriptsize $-1$}}
\put(37,15){\makebox(0,0)[t]{\scriptsize $-2$}}
\put(57,15){\makebox(0,0)[t]{\scriptsize $-3$}}
\put(77,15){\makebox(0,0)[t]{\scriptsize $-4$}}
\put(97,15){\makebox(0,0)[t]{\scriptsize $-5$}}
\put(117,15){\makebox(0,0)[t]{\scriptsize $-6$}}
\put(137,15){\makebox(0,0)[t]{\scriptsize $-7$}}
\put(157,15){\makebox(0,0)[t]{\scriptsize $-8$}}
\put(177,15){\makebox(0,0)[t]{\scriptsize $-9$}}
\put(197,15){\makebox(0,0)[t]{\scriptsize $-10$}}
\put(217,15){\makebox(0,0)[t]{\scriptsize $-11$}}
\put(237,15){\makebox(0,0)[t]{\scriptsize $-12$}}
\put(257,15){\makebox(0,0)[t]{\scriptsize $-13$}}

\put(265,40){\makebox(0,0){$\cdots$}}
\put(265,80){\makebox(0,0){$\cdots$}}

\put(-5,40){\makebox(0,0)[t]{\scriptsize $D_2$}}
\put(-5,85){\makebox(0,0)[t]{\scriptsize $D_1$}}

{\linethickness{0.02pt}
\put(15,60){\line(1,0){265}}}

\thicklines
\put(60,60){\line(-1,1){40}}
\put(80,60){\line(1,1){40}}
\put(100,60){\line(1,1){40}}
\put(120,60){\line(1,1){40}}
\put(180,60){\line(0,1){40}}
\put(200,60){\line(1,1){40}}
\put(220,60){\line(1,1){40}}

\put(30,60){\oval(20,20)[t]}
\put(150,60){\oval(20,20)[t]}
\put(50,100){\oval(20,20)[b]}
\put(90,100){\oval(20,20)[b]}
\put(210,100){\oval(20,20)[b]}

\put(70,20){\oval(20,20)[t]}
\put(110,20){\oval(20,20)[t]}
\put(150,20){\oval(20,20)[t]}
\put(50,60){\oval(20,20)[b]}
\put(90,60){\oval(20,20)[b]}

\put(20,20){\line(0,1){40}}
\put(40,20){\line(2,1){80}}
\put(180,20){\line(-1,1){40}}

\put(200,20){\line(-1,1){40}}
\put(220,20){\line(-1,1){40}}
\put(240,20){\line(-1,1){40}}
\put(260,20){\line(-1,1){40}}

\put(20,20){\line(0,1){40}}
\put(40,20){\line(2,1){80}}
\put(180,20){\line(-1,1){40}}
\put(200,20){\line(-1,1){40}}

\end{picture}
\end{center}

\vspace{-.1cm}

\noindent which yields the diagram

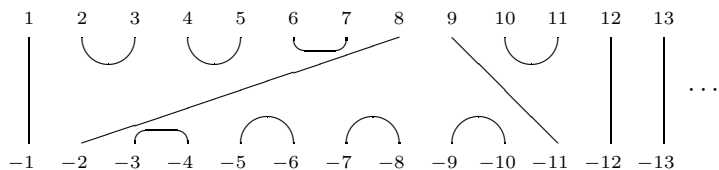
\begin{center}
\begin{picture}(290,60)(0,10)

\put(17,15){\makebox(0,0)[t]{\scriptsize $-1$}}
\put(37,15){\makebox(0,0)[t]{\scriptsize $-2$}}
\put(57,15){\makebox(0,0)[t]{\scriptsize $-3$}}
\put(77,15){\makebox(0,0)[t]{\scriptsize $-4$}}
\put(97,15){\makebox(0,0)[t]{\scriptsize $-5$}}
\put(117,15){\makebox(0,0)[t]{\scriptsize $-6$}}
\put(137,15){\makebox(0,0)[t]{\scriptsize $-7$}}
\put(157,15){\makebox(0,0)[t]{\scriptsize $-8$}}
\put(177,15){\makebox(0,0)[t]{\scriptsize $-9$}}
\put(197,15){\makebox(0,0)[t]{\scriptsize $-10$}}
\put(217,15){\makebox(0,0)[t]{\scriptsize $-11$}}
\put(237,15){\makebox(0,0)[t]{\scriptsize $-12$}}
\put(257,15){\makebox(0,0)[t]{\scriptsize $-13$}}

\put(20,65){\makebox(0,0)[b]{\scriptsize $1$}}
\put(40,65){\makebox(0,0)[b]{\scriptsize $2$}}
\put(60,65){\makebox(0,0)[b]{\scriptsize $3$}}
\put(80,65){\makebox(0,0)[b]{\scriptsize $4$}}
\put(100,65){\makebox(0,0)[b]{\scriptsize $5$}}
\put(120,65){\makebox(0,0)[b]{\scriptsize $6$}}
\put(140,65){\makebox(0,0)[b]{\scriptsize $7$}}
\put(160,65){\makebox(0,0)[b]{\scriptsize $8$}}
\put(180,65){\makebox(0,0)[b]{\scriptsize $9$}}
\put(200,65){\makebox(0,0)[b]{\scriptsize $10$}}
\put(220,65){\makebox(0,0)[b]{\scriptsize $11$}}
\put(240,65){\makebox(0,0)[b]{\scriptsize $12$}}
\put(260,65){\makebox(0,0)[b]{\scriptsize $13$}}

\put(20,20){\line(0,1){40}}
\put(40,20){\line(3,1){120}}
\put(220,20){\line(-1,1){40}}
\put(240,20){\line(0,1){40}}
\put(260,20){\line(0,1){40}}

\put(70,20){\oval(20,10)[t]}
\put(110,20){\oval(20,20)[t]}
\put(150,20){\oval(20,20)[t]}
\put(190,20){\oval(20,20)[t]}

\put(50,60){\oval(20,20)[b]}
\put(90,60){\oval(20,20)[b]}
\put(130,60){\oval(20,10)[b]}
\put(210,60){\oval(20,20)[b]}

\put(275,40){\makebox(0,0){$\cdots$}}

\end{picture}
\end{center}

\vspace{-.1cm}

To ascertain that friezes are closed under composition, that this
composition is associative, and that the unit frieze is indeed a
unit with respect to composition, i.e.\ to ascertain that friezes
make a monoid, we need a more formal definition of composition.
Formally, we may define composition of friezes either in a
geometrical style (see \cite{DP01}), or in a set-theoretical
style, as a peculiar composition of equivalence relations (see
\cite{DP021} or \cite{DP022}). We don't have space here to go into
these formal matters, which have already been treated elsewhere.
So we take for granted that friezes make a monoid.

An important aspect of this matter is that
no circles, or closed loops, made of caps from one frieze and
cups from the other can arise in composition of friezes.
This particular fact is ascertained easily from the
distribution of even and odd numbers in the end points of cups
and caps.

The adjunctional friezes of the present paper are obtained from the
friezes of \cite{DP01} by permitting infinitely many cups and caps,
by forgetting about circles and by requiring that
cups be of the form $[2k\pl 2, 2k\pl 3]$ and
caps of the form $[-(2k\pl 2),-(2k\pl 1)]$.
Without this condition on cups and caps, circles may arise
in composition, and according to how we treat them we
obtain in \cite{DP01} various kinds of equivalences of friezes,
and various kinds of monoids.
The strictest, $\cal L$ kind, records in what regions of the diagram
the circles are located.
The $\cal K$ kind just counts the number of circles, and the loosest,
$\cal J$ kind, ignores circles. Since we don't have
circles in the adjunctional friezes of this paper,
the various notions of equivalence of friezes of \cite{DP01} will
coincide for them. The friezes of \cite{DP01} are tied to
self-adjoint situations, where an endofunctor is adjoint to
itself, while the friezes introduced here are tied to
arbitrary adjoint situations.

Something analogous to adjunctional and other friezes may be found in
\cite{B37}, \cite{EK66}, \cite{KM71}, \cite{Y88}, \cite{W88},
\cite{T89}, and in many other papers in knot theory following Jones'
approach to knot and link invariants, which we mentioned in the introduction
and at the beginning of this section.

Still another possibility to define composition of friezes is to
rely on the isomorphism with \Ss\ (see the end
of this section). We are now going to establish this isomorphism.

A transversal segment of a frieze is {\it odd} when
its end points are two odd integers, and analogously with ``odd''
replaced by ``even''. Every transversal segment of a frieze
is either odd or even.

Let the {\it successor} of a transversal segment in a frieze
be the next transversal segment on the right-hand side in the diagram.
For example, in the frieze $D_1$, which we
had as our first example above,
the successor of $[-3,1]$ is $[-4,6]$. The successor of an
odd transversal segment $[-(2k_1\pl 1), 2k_2\pl 1]$ is an even
transversal segment $[-(2k_1\pl 2), 2k_2\pl 2k_3 \pl 2]$, and
the successor of an even
transversal segment $[-2k_1, 2k_2]$ is an odd
transversal segment $[-(2k_1\pl 2k_3 \pl 1), 2k_2\pl 1]$.

Let a {\it transversal pair} in a frieze
be an odd transversal segment
$[-(2k_1\pl 1), 2k_2\pl 1]$ and its
successor $[-(2k_1\pl 2), 2k_2\pl 2k_3 \pl 2]$. Let us say that
$n\in \N$ is {\it covered} by this transversal pair when
$k_2\leq n\leq k_2\pl k_3$, and let us say that this
transversal pair {\it assigns} $k_1$ to $n$. In the example with
$D_1$ above, the transversal pair made of $[-3,1]$ and $[-4,6]$
covers $0$, $1$ and $2$, and it assigns $1$ to these three numbers.
In every frieze, every $n\in \N$ is covered by
exactly one transversal pair, and this transversal pair
assigns to every $n$ a single number $k_1$. What happens should
be clear from the following adaptation of the diagram of $D_1$:

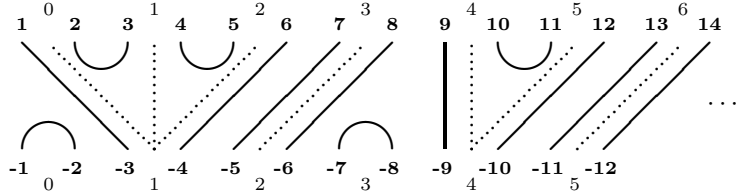
\vspace{-.4cm}

\begin{center}
\begin{picture}(290,75)(0,10)

\put(19,15){\makebox(0,0)[t]{\scriptsize \bf -1}}
\put(30,9){\makebox(0,0)[t]{\scriptsize $0$}}
\put(39,15){\makebox(0,0)[t]{\scriptsize \bf -2}}
\put(59,15){\makebox(0,0)[t]{\scriptsize \bf -3}}
\put(70,9){\makebox(0,0)[t]{\scriptsize $1$}}
\put(79,15){\makebox(0,0)[t]{\scriptsize \bf -4}}
\put(99,15){\makebox(0,0)[t]{\scriptsize \bf -5}}
\put(110,9){\makebox(0,0)[t]{\scriptsize $2$}}
\put(119,15){\makebox(0,0)[t]{\scriptsize \bf -6}}
\put(139,15){\makebox(0,0)[t]{\scriptsize \bf -7}}
\put(150,9){\makebox(0,0)[t]{\scriptsize $3$}}
\put(159,15){\makebox(0,0)[t]{\scriptsize \bf -8}}
\put(179,15){\makebox(0,0)[t]{\scriptsize \bf -9}}
\put(190,9){\makebox(0,0)[t]{\scriptsize $4$}}
\put(199,15){\makebox(0,0)[t]{\scriptsize \bf -10}}
\put(219,15){\makebox(0,0)[t]{\scriptsize \bf -11}}
\put(229,9){\makebox(0,0)[t]{\scriptsize $5$}}
\put(239,15){\makebox(0,0)[t]{\scriptsize \bf -12}}

\put(20,65){\makebox(0,0)[b]{\scriptsize \bf 1}}
\put(30,75){\makebox(0,0)[t]{\scriptsize $0$}}
\put(40,65){\makebox(0,0)[b]{\scriptsize \bf 2}}
\put(60,65){\makebox(0,0)[b]{\scriptsize \bf 3}}
\put(70,75){\makebox(0,0)[t]{\scriptsize $1$}}
\put(80,65){\makebox(0,0)[b]{\scriptsize \bf 4}}
\put(100,65){\makebox(0,0)[b]{\scriptsize \bf 5}}
\put(110,75){\makebox(0,0)[t]{\scriptsize $2$}}
\put(120,65){\makebox(0,0)[b]{\scriptsize \bf 6}}
\put(140,65){\makebox(0,0)[b]{\scriptsize \bf 7}}
\put(150,75){\makebox(0,0)[t]{\scriptsize $3$}}
\put(160,65){\makebox(0,0)[b]{\scriptsize \bf 8}}
\put(180,65){\makebox(0,0)[b]{\scriptsize \bf 9}}
\put(190,75){\makebox(0,0)[t]{\scriptsize $4$}}
\put(200,65){\makebox(0,0)[b]{\scriptsize \bf 10}}
\put(220,65){\makebox(0,0)[b]{\scriptsize \bf 11}}
\put(230,75){\makebox(0,0)[t]{\scriptsize $5$}}
\put(240,65){\makebox(0,0)[b]{\scriptsize \bf 12}}
\put(260,65){\makebox(0,0)[b]{\scriptsize \bf 13}}
\put(270,75){\makebox(0,0)[t]{\scriptsize $6$}}
\put(280,65){\makebox(0,0)[b]{\scriptsize \bf 14}}

\put(285,37){\makebox(0,0){$\cdots$}}

\multiput(70,20)(-2,2){20}{\circle*{.1}}

\multiput(70,20)(0,3){14}{\circle*{.1}}

\multiput(70,20)(2,2){20}{\circle*{.1}}

\multiput(110,20)(2,2){20}{\circle*{.1}}

\multiput(190,20)(0,3){14}{\circle*{.1}}

\multiput(190,20)(2,2){20}{\circle*{.1}}

\multiput(230,20)(2,2){20}{\circle*{.1}}

\thicklines
\put(60,20){\line(-1,1){40}}
\put(80,20){\line(1,1){40}} \put(100,20){\line(1,1){40}}
\put(120,20){\line(1,1){40}} \put(180,20){\line(0,1){40}}
\put(200,20){\line(1,1){40}} \put(220,20){\line(1,1){40}}
\put(240,20){\line(1,1){40}}

\put(30,20){\oval(20,20)[t]}
\put(150,20){\oval(20,20)[t]}
\put(50,60){\oval(20,20)[b]}
\put(90,60){\oval(20,20)[b]}
\put(210,60){\oval(20,20)[b]}

\end{picture}
\end{center}

\vspace{.1cm}

\noindent So for every frieze $D$ we can define an
order-preserving function $\varphi(D)\!:\N \str \N$ by mapping $n$
to the number assigned to $n$ by the transversal pair of $D$
covering $n$.

Conversely, for every order-preserving function $\varphi\!:\N\str
\N$ we can define a frieze $D(\varphi)$ whose transversal segments
are obtained as follows. For every $k$ such that there is an $n$
for which $\varphi(n)=k$, we have the transversal segments
$[-(2k\pl 1), 2\min\{n\mid \varphi(n)=k\}\pl 1]$ and $[-(2k\pl 2),
2\max\{n\mid \varphi(n)=k\}\pl 2]$; these two segments make a
transversal pair. It is easy to check that $D(\varphi(D))=D$ and
$\varphi(D(\varphi))=\varphi$, so that we have a bijection between
friezes and order-preserving endomorphisms of \N.

It is clear that for the unit frieze \mj\ we have that
$\varphi(\mj)$ is the identity function on \N. We
also have the following lemma.

\vspace{.3cm}

\noindent $\varphi$ {\sc Lemma}. $\varphi(D_2\cirk D_1)=\varphi(D_2)
\cirk \varphi(D_1)$.

\vspace{.2cm}

\noindent {\it Proof.} Let $n\in \N$ be covered by the transversal pair
of $D_1$ whose segments are
$[-(2k_4\pl 2k_5\pl 1), 2k_2\pl 1]$
and $[-(2k_4\pl 2k_5\pl 2), 2k_2\pl 2k_3 \pl 2]$, and
let $k_4\pl k_5$ be covered by the transversal pair of
$D_2$ whose segments are
$[-(2k_1\pl 1), 2k_4\pl 1]$
and $[-(2k_1\pl 2), 2k_2\pl 2k_4\pl 2k_5\pl 2k_6\pl 2]$, so that
$\varphi(D_2)(\varphi(D_1)(n))=k_1$. Then $n$ is covered by the
transversal pair of $D_2\cirk D_1$ whose segments are
$[-(2k_1\pl 1), 2k_2\mn 2k_7\pl 1]$
and $[-(2k_1\pl 2), 2k_2\pl 2k_3\pl 2k_8 \pl 2]$, and
so $\varphi(D_2\cirk D_1)(n)=k_1$. \qed

\vspace{.3cm}

Hence we have that $\varphi$ establishes an isomorphism between
the monoid of friezes and the monoid \Ss\ of order-preserving
endomorphisms of \N.

We could have defined formally
composition of friezes by relying on the bijectivity of $\varphi$.
The composition $D_2\cirk D_1$ of the friezes $D_1$ and $D_2$
could be defined as $D(\varphi(D_2)\cirk \varphi(D_1))$. With this
definition it is trivial to establish that friezes make a monoid
isomorphic to \Ss. But then it remains to establish that composition so
defined is the same notion we find in other possible formal definitions of
composition, which we mentioned above.

The isomorphism of the monoid of friezes with the monoid \Ss\ of
order-preserving endomorphisms of \N, which was established in
this section, can be relied upon in investigating an `augmented'
simplicial category, whose objects are the finite ordinals
together with \N\ (i.e.\ the ordinal $\omega$), and whose arrows
are the order-preserving functions. This point was raised by an
anonymous referee of this paper, who asked whether there is a
faithful functor from the augmented simplicial category to a
category with the same objects whose arrows are instances of an
appropriately modified notion of frieze. It seems likely that the
question can be answered positively, but to give a precise answer
would lead us too far afield, and we leave the matter for future
research.

\section{Adjunction and friezes}

In this section we show how friezes are tied to the free adjunction
generated by $(0,1)$, which we introduced in Section 5. The connection
between friezes and adjunction is made, more or less implicitly, in
\cite{KM71}, \cite{Y88}, \cite{FY89}, \cite{FY92} and \cite{KD99}.

For $n,m\in \N$, a frieze is said to be of {\it type} $(n,m)$ when for
every $k\in \N^+=\N\mn \{0\}$ we have a transversal segment
$[-(m\pl k), n\pl k]$ in this frieze.
For example, the frieze $D_1$ from the preceding section is of
type $(11,9)$.
Not all friezes have a type, and when
they have one, they are said to be of {\it finite type}.
Note that types of friezes are not unique;
a frieze of type $(n,m)$ is also of type $(n\pl k,m\pl k)$.
It is clear that if $D_1$ is a frieze of type $(n,m)$ and $D_2$ a
frieze of type $(m,l)$, then $D_2\cirk D_1$ is a frieze of type
$(n,l)$.

Let \Da\ be the category whose objects are all odd natural numbers,
and whose arrows between $n$ and $m$ are all friezes of type
$(n,m)$ indexed by $(n,m)$. We index these friezes by $(n,m)$ to
ensure that every arrow has a single source $n$ and a single target
$m$ (as we said above, every frieze of finite type has infinitely
many different types). The category \Db\ is defined analogously save that
its objects are all even numbers (including $0$). We will show that
these categories are isomorphic respectively to the categories
\A\ and \B\ of the free adjunction generated by $(0,1)$.

We define a functor \Ea\ from \A\ to \Da\ and a functor \Eb\ from
\B\ to \Db. On objects we have

\[
\begin{array}{l}

\Ea(F(GF)^n\O)=2n\pl 1,
\\[.05cm]
\Eb((GF)^n\O)=2n.

\end{array}
\]

\noindent Next we define inductively \Ea\ and \Eb\
on the arrow terms of \A\ and \B:

\begin{description}

\item[\hspace{.4cm}] $\Ea(\mj_a)$ is the unit frieze
\mj\ indexed by $(\Ea(a),\Ea(a))$;
\vspace{-.25cm}
\item[\hspace{.4cm}] $\Eb(\mj_b)$ is the unit frieze
\mj\ indexed by $(\Eb(b),\Eb(b))$;
\vspace{-.05cm}
\item[\hspace{.4cm}] $\Ea(f_2\cirk f_1)$ is $\Ea(f_2)\cirk \Ea(f_1)$;
\vspace{-.25cm}
\item[\hspace{.4cm}] $\Eb(g_2\cirk g_1)$ is $\Eb(g_2)\cirk \Eb(g_1)$;
\vspace{-.05cm}
\item[\hspace{.4cm}] $\Ea(Fg)$ is $\Eb(g)$ with its index $(n,m)$
replaced by $(n\pl 1, m\pl 1)$;
\vspace{-.25cm}
\item[\hspace{.4cm}] $\Eb(Gf)$ is $\Ea(f)$ with its index $(n,m)$
replaced by $(n\pl 1, m\pl 1)$;
\vspace{-.05cm}
\item[\hspace{.4cm}] if $\Ea(f)$ is indexed by $(n,m)$, then
$\Ea(\fia f)$ is the frieze indexed by $(n\pl 2,m)$ obtained from
the frieze $\Ea(f)$ by replacing all the
transversal segments $[-(m\pl k), n\pl k]$, for every $k\in \N^+$,
by the cup $[n\pl 1, n\pl 2]$ and by the transversal segments
$[-(m\pl k), n\pl k\pl 2]$;
\vspace{-.25cm}
\item[\hspace{.4cm}] if $\Eb(g)$ is indexed by $(n,m)$, then
$\Eb(\gac g)$ is the frieze indexed by $(n,m\pl 2)$ obtained from
the frieze $\Eb(g)$ by replacing all the
transversal segments $[-(m\pl k), n\pl k]$, for every $k\in \N^+$,
by the cap $[-(m\pl 2),-(m\pl 1)]$ and by the transversal segments
$[-(m\pl k\pl 2), n\pl k]$.

\end{description}

\noindent Let us illustrate the last two clauses, for \fia\ and \gac:

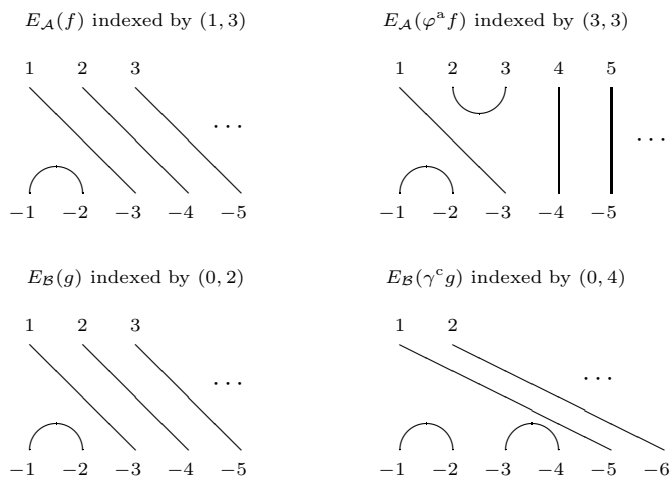
\begin{center}
\begin{picture}(220,90)

\put(-3,15){\makebox(0,0)[t]{\scriptsize $-1$}}
\put(17,15){\makebox(0,0)[t]{\scriptsize $-2$}}
\put(37,15){\makebox(0,0)[t]{\scriptsize $-3$}}
\put(57,15){\makebox(0,0)[t]{\scriptsize $-4$}}
\put(77,15){\makebox(0,0)[t]{\scriptsize $-5$}}

\put(0,65){\makebox(0,0)[b]{\scriptsize $1$}}
\put(20,65){\makebox(0,0)[b]{\scriptsize $2$}}
\put(40,65){\makebox(0,0)[b]{\scriptsize $3$}}

\put(75,45){\makebox(0,0){$\cdots$}}

\put(10,20){\oval(20,20)[t]}

\put(40,20){\line(-1,1){40}}
\put(60,20){\line(-1,1){40}}
\put(80,20){\line(-1,1){40}}

\put(40,85){\makebox(0,0){\scriptsize{$\Ea(f)$ indexed by $(1,3)$}}}

\put(137,15){\makebox(0,0)[t]{\scriptsize $-1$}}
\put(157,15){\makebox(0,0)[t]{\scriptsize $-2$}}
\put(177,15){\makebox(0,0)[t]{\scriptsize $-3$}}
\put(197,15){\makebox(0,0)[t]{\scriptsize $-4$}}
\put(217,15){\makebox(0,0)[t]{\scriptsize $-5$}}

\put(140,65){\makebox(0,0)[b]{\scriptsize $1$}}
\put(160,65){\makebox(0,0)[b]{\scriptsize $2$}}
\put(180,65){\makebox(0,0)[b]{\scriptsize $3$}}
\put(200,65){\makebox(0,0)[b]{\scriptsize $4$}}
\put(220,65){\makebox(0,0)[b]{\scriptsize $5$}}

\put(235,40){\makebox(0,0){$\cdots$}}

\put(150,20){\oval(20,20)[t]}
\put(170,60){\oval(20,20)[b]}

\put(180,20){\line(-1,1){40}}
\put(200,20){\line(0,1){40}}
\put(220,20){\line(0,1){40}}

\put(180,85){\makebox(0,0){\scriptsize{$\Ea(\fia f)$ indexed by $(3,3)$}}}

\end{picture}
\end{center}

\vspace{-.3cm}

\begin{center}
\begin{picture}(220,90)

\put(-3,15){\makebox(0,0)[t]{\scriptsize $-1$}}
\put(17,15){\makebox(0,0)[t]{\scriptsize $-2$}}
\put(37,15){\makebox(0,0)[t]{\scriptsize $-3$}}
\put(57,15){\makebox(0,0)[t]{\scriptsize $-4$}}
\put(77,15){\makebox(0,0)[t]{\scriptsize $-5$}}

\put(0,65){\makebox(0,0)[b]{\scriptsize $1$}}
\put(20,65){\makebox(0,0)[b]{\scriptsize $2$}}
\put(40,65){\makebox(0,0)[b]{\scriptsize $3$}}

\put(75,45){\makebox(0,0){$\cdots$}}

\put(10,20){\oval(20,20)[t]}

\put(40,20){\line(-1,1){40}}
\put(60,20){\line(-1,1){40}}
\put(80,20){\line(-1,1){40}}

\put(40,85){\makebox(0,0){\scriptsize{$\Eb(g)$ indexed by $(0,2)$}}}

\put(137,15){\makebox(0,0)[t]{\scriptsize $-1$}}
\put(157,15){\makebox(0,0)[t]{\scriptsize $-2$}}
\put(177,15){\makebox(0,0)[t]{\scriptsize $-3$}}
\put(197,15){\makebox(0,0)[t]{\scriptsize $-4$}}
\put(217,15){\makebox(0,0)[t]{\scriptsize $-5$}}
\put(237,15){\makebox(0,0)[t]{\scriptsize $-6$}}

\put(140,65){\makebox(0,0)[b]{\scriptsize $1$}}
\put(160,65){\makebox(0,0)[b]{\scriptsize $2$}}

\put(215,47){\makebox(0,0){$\cdots$}}

\put(150,20){\oval(20,20)[t]}
\put(190,20){\oval(20,20)[t]}

\put(220,20){\line(-2,1){80}}
\put(240,20){\line(-2,1){80}}

\put(180,85){\makebox(0,0){\scriptsize{$\Eb(\gac g)$ indexed by $(0,4)$}}}

\end{picture}
\end{center}

\vspace{-.3cm}

We can check by induction on the length of derivation that if $f_1=f_2$ in
\A, then $\Ea(f_1)=\Ea(f_2)$ in \Da, and that if $g_1=g_2$ in \B, then
$\Eb(g_1)=\Eb(g_2)$ in \Db. For example, for the first $(\fia \gac)$
equation we have

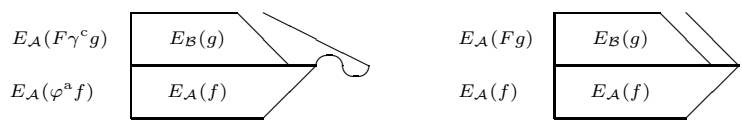
\vspace{-.1cm}

\begin{center}
\begin{picture}(260,60)

\put(30,10){\line(1,0){50}}
\put(80,10){\line(1,1){20}}
\put(30,30){\line(1,0){70}}
\put(30,50){\line(1,0){40}}
\put(90,30){\line(-1,1){20}}
\put(120,30){\line(-2,1){40}}
\put(30,10){\line(0,1){40}}

\put(105,30){\oval(10,8)[t]}
\put(115,30){\oval(10,8)[b]}

\put(0,20){\makebox(0,0){\scriptsize{$\Ea(\fia f)$}}}
\put(3,40){\makebox(0,0){\scriptsize{$\Ea(F\gac g)$}}}

\put(55,20){\makebox(0,0){\scriptsize{$\Ea(f)$}}}
\put(55,40){\makebox(0,0){\scriptsize{$\Eb(g)$}}}

\put(190,10){\line(1,0){50}}
\put(240,10){\line(1,1){20}}
\put(190,30){\line(1,0){70}}
\put(190,50){\line(1,0){40}}
\put(250,30){\line(-1,1){20}}
\put(190,10){\line(0,1){40}}

\put(260,30){\line(-1,1){20}}

\put(165,20){\makebox(0,0){\scriptsize{$\Ea(f)$}}}
\put(168,40){\makebox(0,0){\scriptsize{$\Ea(Fg)$}}}

\put(215,20){\makebox(0,0){\scriptsize{$\Ea(f)$}}}
\put(215,40){\makebox(0,0){\scriptsize{$\Eb(g)$}}}

\end{picture}
\end{center}

\vspace{-.2cm}

\noindent This shows that the triangular equations of adjunctions,
the essential equations of adjunctions, to which the equations
$(\fia \gac)$ correspond, are about ``straightening a sinuosity''
(cf.\ \cite{KD99}, \S 4.10.1), and this straightening is based on
planar ambient isotopies of knot theory (cf.\ \cite{BZ85}, \S
1.A).

Since it is clear that \Ea\ and \Eb\ preserve identities and
composition, with \Ea\ and \Eb\ we have indeed functors from \A\
to \Da\ and from \B\ to \Db\ respectively. We will establish that
these functors are isomorphisms. First we prove the following
lemma.

\vspace{.3cm}

\noindent {\sc Surjectivity Lemma}. {\it The functors \Ea\ and \Eb\
are onto on arrows.}

\vspace{.2cm}

\noindent {\it Proof.} Take an arbitrary arrow $D$ of \Da\ or \Db\ of
type $(n,m)$. We construct a term $h$ of \A\ or \B\ such that
$\Ea(h)=D$ or $\Eb(h)=D$, as appropriate, by induction on $n\pl m$.

Suppose $n\pl m=0$; so $n=m=0$. Then $D$ must be the unit frieze
\mj\ indexed by (0,0), which is an arrow of \Db, and $h$ is the
term $\mj_{\O}\!:\O\str \O$ of \B.

Suppose $n+m>0$, and suppose $n$ and $m$ are even. Then $D$ is an arrow
of \Db, and in $D$ we either have a transversal segment $[-m,n]$, or
a cap $[-m,-(m\mn 1)]$. In the first case, by the induction hypothesis,
we have constructed a term $h'$ of \A\ such that $\Ea(h')$ is the
arrow $D$ of \Da\ of type $(n\mn 1, m\mn 1)$. Then $h$ is $Gh'$. In the
second case, by the induction hypothesis, we have constructed the term
$h'$ of \B\ such that $\Eb(h')$ is the arrow $D'$ of \Db\ of type
$(n,m\mn 2)$, and $D$ is obtained from $D'$ by replacing the
all the transversal segments $[-(m\mn2 \pl k), n\pl k]$, for every
$k\in \N^+$, by the cap $[-m,-(m\mn 1)]$ and by the transversal segments
$[-(m\pl k), n\pl k]$. Then $h$ is $\gac h'$. We proceed analogously,
in a dual manner, when $n$ and $m$ are odd. \qed

\vspace{.3cm}

For the next lemma we rely on the normal form of terms
of \A\ and \B\ defined in Section 5.

\vspace{.3cm}

\noindent {\sc Auxiliary Lemma.} (\A) {\it If $f_1,f_2\!:a_1\str
a_2$ are terms of \A\ in normal form and $\Ea(f_1)=\Ea(f_2)$ in
\Da, then $f_1$ is the same term as $f_2$.}

(\B) {\it If $g_1,g_2\!:b_1\str b_2$ are terms of \B\ in normal
form and $\Eb(g_1)=\Eb(g_2)$ in \Db, then $g_1$ is the same term
as $g_2$.}

\vspace{.2cm}

\noindent {\it Proof.} (\A) Suppose $f_1$, which
is $X^n\ldots X^1\mj_{\O}$, and
$f_2$, which is $Y^m\ldots Y^1\mj_{\O}$, are different terms of
\A\ of the same type $a_1\str a_2$. Suppose for some $i$ such that
$1\leq i\leq n$ and $1\leq i\leq m$ we have that $X_i$ is different from
$Y_i$. Let $j$ be the least such $i$. This means either that one of
$X_j$ and $Y_j$ is $F$, while the other is $\gac$, or that one of
$X_j$ and $Y_j$ is $G$, while the other is $\fia$. In both cases we
obtain that $\Ea(f_1)$ cannot be equal to $\Ea(f_2)$. If for every $i$
such that $1\leq i\leq n$ we have that $X_i$ is identical to $Y_i$,
then $f_1$ and $f_2$ can differ only if $n<m$.
But then $f_1$ and $f_2$ cannot be of the same type.
The proof of (\B) is analogous. \qed

\vspace{.3cm}

As a consequence of the last lemma, and of the functoriality of
\Ea\ and \Eb, we obtain that the normal form of terms of \A\ and \B\
is unique. We reason as in Sections 3 and 4. We also obtain the
following lemma, by reasoning as in the proofs of the Injectivity
Lemmata of Sections 3 and 4.

\vspace{.3cm}

\noindent {\sc Injectivity Lemma}. {\it If $f_1,f_2\!:a_1\str a_2$
are terms of \A\ and $\Ea(f_1)=\Ea(f_2)$ in \Da, then $f_1=f_2$ in
\A. If $g_1,g_2\!:b_1\str b_2$ are terms of \B\ and
$\Eb(g_1)=\Eb(g_2)$ in \Db, then $g_1=g_2$ in \B.}

\vspace{.3cm}

Since \Ea\ and \Eb\ are bijections on objects, this lemma says that
these functors are one-one on arrows. So, together with the Surjectivity
Lemma, this yields that they are isomorphisms.

If the free adjunction is generated by $(1,0)$ instead of $(0,1)$,
then in \A\ we have the objects $(FG)^n\O$ and in \B\ the objects
$G(FG)^n\O$, so that in \Da\ the objects are even and in \Db\ odd
natural numbers. The definition of frieze would then have
cups $[2k\pl 1, 2k\pl 2]$ and caps $[-(2k\pl 3), -(2k\pl 2)]$.
These new friezes correspond to relations converse to
order-preserving functions. We could again establish
isomorphisms between \A\ and \Da\ and between \B\ and \Db,
analogous to those above.
(This is done in \cite{KD99}, \S 4.10,
but in a manner somewhat less detailed than here, in particular
with respect to friezes; the idea is, however, the same.)
If the free adjunction is generated by
$(1,1)$, then the category \A\ is the disjoint union of the categories
\A\ generated by $(0,1)$ and $(1,0)$, and analogously for \B.
We have then again isomorphisms with the appropriately
defined categories \Da\ and \Db.

\section{Adjunction and \Sc}

We will now show that the category \Db\ of the preceding section
is isomorphic to the simplicial category \Sc. We define first a
functor $S$ from \Db\ to \Sc. On objects, we have that $S(2n)=n$.
On arrows, we have that if $D$ is a frieze indexed by $(2n,2m)$,
then $S(D)$ is the order-preserving function from $n$ to $m$
obtained by restricting the domain of the order-preserving
endomorphism $\varphi(D)\!:\N\str \N$ (see Section 6) to $n$, and
the codomain to $m$.

In a frieze $D$ of type $(2n,2m)$, for some $k<2m$ we have
a transversal segment $[-k,1]$, and, provided $n>0$, for some
$l\leq 2m$ we have a transversal segment $[-l,2n]$. Let us call the
transversal segments in between these two transversal segments,
including these two segments, the {\it specific} transversal
segments of $D$. If $n>0$, then the number of specific
transversal segments is an even number greater than or equal to $2$,
and if $n=0$, then $2n$ is not the end point of a transversal
segment, and we have zero specific transversal segments. Specific
transversal segments are those out of which we take the transversal
pairs that determine the value of $S(D)$ for every number in
$\{0,\ldots ,n\mn 1\}$. This value is a number in
$\{0,\ldots ,m\mn 1\}$. If $n=0$, then there are no transversal
pairs made of specific transversal segments,
and $S(D)$ is the empty function.

For the unit frieze \mj\ we have that $S(\mj)$ is the identity
function, and
we check, as in the $\varphi$ Lemma of Section 6, that $S$ preserves
composition of friezes. So $S$ is indeed a functor from \Db\ to \Sc.

To show that $S$ is an isomorphism, we define a functor $D$ from
\Sc\ to \Db. On objects, we have that $D(n)=2n$, and on arrows
$\varphi\!: n\str m$ we define $D(\varphi)$ of type $(2n,2m)$ as
in Section 6. According to what we know from Section 6, it is
clear that these two functors establish an isomorphism between
\Db\ and \Sc.

Note that the category \Da\ of the preceding section is isomorphic
to the subcategory of \Db\ for whose arrows $D$ of type $(2n,2m)$
we have in $D$ the transversal segment $[-2m,2n]$. Since this
segment is paired with $[-(2m\mn 1),k]$, for some $k<2n$, to make
a transversal pair, the category \Da\ is isomorphic to the
subcategory of \Sc\ whose order-preserving functions are
last-element-preserving, i.e.\ the order-preserving function
$\varphi\!: n\str m$ maps $n\mn 1$ to $m\mn 1$.

Since \A\ and \B\ are isomorphic to \Da\ and \Db\ respectively, we
have that \A\ is isomorphic to a subcategory of \B. This is the
subcategory of \B\ whose objects are of the form $(GF)^n\O$ for
some $n\geq 1$, and whose arrows are arrows of \B\ of the form
$Gf$ (cf.\ \cite{KD99}, \S 5.2.2).

The isomorphism of \Db\ with \Sc\ yields the isomorphism of \B\
with \Sc. The existence of this last isomorphism may be extracted
from \cite{A74} (Corollary 2.8; see also \cite{SS86}).

\section{\Sn\ and Temperley-Lieb algebras}

We can now establish that the monoids \Sn\ are submonoids of the
monoids tied to Temperley-Lieb algebras.

Let \Bn\ be the monoid of endomorphisms $g\!:(GF)^n\O \str
(GF)^n\O$ of \B\ (this monoid should not be confused with a braid
group, which is often named similarly), and let \Dbn\ be the
monoid of endomorphisms $D$ of type $(2n,2n)$ of \Db. (It is clear
that such endomorphisms are closed under composition.) According
to what we have established in Section 7, the monoids \Bn\ and
\Dbn\ are isomorphic, and they are both isomorphic to the monoid
\Sn. Here the right-forking term $p^i$ of \Sn\ corresponds to the
frieze

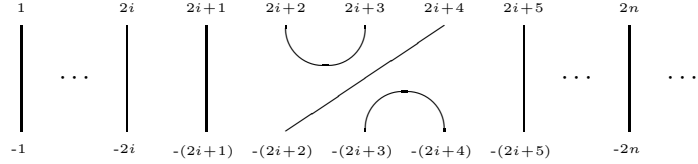
\begin{center}
\begin{picture}(280,60)(0,10)

\put(19,15){\makebox(0,0)[t]{\tiny -$1$}}
\put(59,15){\makebox(0,0)[t]{\tiny -$2i$}}
\put(89,15){\makebox(0,0)[t]{\tiny -$(2i\pl 1)$}}
\put(119,15){\makebox(0,0)[t]{\tiny -$(2i\pl 2)$}}
\put(149,15){\makebox(0,0)[t]{\tiny -$(2i\pl 3)$}}
\put(179,15){\makebox(0,0)[t]{\tiny -$(2i\pl 4)$}}
\put(209,15){\makebox(0,0)[t]{\tiny -$(2i\pl 5)$}}
\put(249,15){\makebox(0,0)[t]{\tiny -$2n$}}

\put(20,65){\makebox(0,0)[b]{\tiny $1$}}
\put(60,65){\makebox(0,0)[b]{\tiny $2i$}}
\put(90,65){\makebox(0,0)[b]{\tiny $2i\pl 1$}}
\put(120,65){\makebox(0,0)[b]{\tiny $2i\pl 2$}}
\put(150,65){\makebox(0,0)[b]{\tiny $2i\pl 3$}}
\put(180,65){\makebox(0,0)[b]{\tiny $2i\pl 4$}}
\put(210,65){\makebox(0,0)[b]{\tiny $2i\pl 5$}}
\put(250,65){\makebox(0,0)[b]{\tiny $2n$}}

\put(20,20){\line(0,1){40}}
\put(60,20){\line(0,1){40}}
\put(90,20){\line(0,1){40}}
\put(210,20){\line(0,1){40}}
\put(250,20){\line(0,1){40}}
\put(120,20){\line(3,2){60}}

\put(165,20){\oval(30,30)[t]}
\put(135,60){\oval(30,30)[b]}

\put(40,40){\makebox(0,0){$\cdots$}}
\put(230,40){\makebox(0,0){$\cdots$}}
\put(270,40){\makebox(0,0){$\cdots$}}

\end{picture}
\end{center}

\vspace{-.1cm}

\noindent which corresponds to the following arrow of \B:

\[
(GF)^{n-i-2}\gac G\fia F(GF)^i\mj_{\O}\!:(GF)^n\O \str (GF)^n\O,
\]

\noindent and the left-forking term $q^i$ corresponds to the frieze

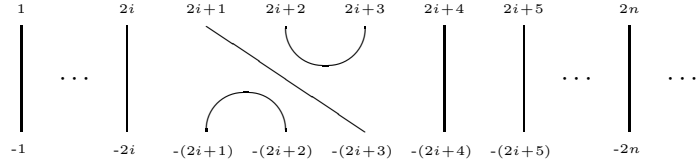
\begin{center}
\begin{picture}(280,60)(0,10)

\put(19,15){\makebox(0,0)[t]{\tiny -$1$}}
\put(59,15){\makebox(0,0)[t]{\tiny -$2i$}}
\put(89,15){\makebox(0,0)[t]{\tiny -$(2i\pl 1)$}}
\put(119,15){\makebox(0,0)[t]{\tiny -$(2i\pl 2)$}}
\put(149,15){\makebox(0,0)[t]{\tiny -$(2i\pl 3)$}}
\put(179,15){\makebox(0,0)[t]{\tiny -$(2i\pl 4)$}}
\put(209,15){\makebox(0,0)[t]{\tiny -$(2i\pl 5)$}}
\put(249,15){\makebox(0,0)[t]{\tiny -$2n$}}

\put(20,65){\makebox(0,0)[b]{\tiny $1$}}
\put(60,65){\makebox(0,0)[b]{\tiny $2i$}}
\put(90,65){\makebox(0,0)[b]{\tiny $2i\pl 1$}}
\put(120,65){\makebox(0,0)[b]{\tiny $2i\pl 2$}}
\put(150,65){\makebox(0,0)[b]{\tiny $2i\pl 3$}}
\put(180,65){\makebox(0,0)[b]{\tiny $2i\pl 4$}}
\put(210,65){\makebox(0,0)[b]{\tiny $2i\pl 5$}}
\put(250,65){\makebox(0,0)[b]{\tiny $2n$}}

\put(20,20){\line(0,1){40}}
\put(60,20){\line(0,1){40}}
\put(180,20){\line(0,1){40}}
\put(210,20){\line(0,1){40}}
\put(250,20){\line(0,1){40}}
\put(150,20){\line(-3,2){60}}

\put(105,20){\oval(30,30)[t]}
\put(135,60){\oval(30,30)[b]}

\put(40,40){\makebox(0,0){$\cdots$}}
\put(230,40){\makebox(0,0){$\cdots$}}
\put(270,40){\makebox(0,0){$\cdots$}}

\end{picture}
\end{center}

\vspace{-.1cm}

\noindent which corresponds to the following arrow of \B:

\[
(GF)^{n-i-2}G\fia F\gac (GF)^i\mj_{\O}\!:(GF)^n\O \str (GF)^n\O.
\]

The monoid \Kn\ has for
every $i \in \{1,\ldots, n\mn 1\}$ a generator $h^i$,
and also the generator $c$.
The terms of \Kn\
are obtained from these generators and \mj\ by closing
under the binary operation \cirk.
The following equations are assumed for \Kn\ besides the monoid
equations $(1)$ and $(2)$ of Section 2:

\[
\begin{array}{ll}
\makebox[5em][l]{$(h1)$} & \makebox[15em][l]
{$h^i \cirk h^j=h^j \cirk h^i, \quad {\mbox{\rm for }} j\pl 1<i$,}\\[.05cm]
(h2) & h^i \cirk h^{i\pm 1} \cirk h^i=h^i,
\end{array}
\]
\vspace{-.2cm}
\[
\begin{array}{ll}
\makebox[5em][l]{$(hc1)$} & \makebox[15em][l]
{$h^i \cirk c=c \cirk h^i$,}\\[.05cm]
(hc2) & h^i \cirk h^i=c \cirk h^i.
\end{array}
\]

\vspace{.2cm}

\noindent The equations $(h1)$, $(h2)$ and $(hc2)$,
which may be derived from Jones' paper \cite{J83} (p. 13),
and which appear in the form above in many works of Kauffman (see
\cite{KL}, \cite{KA97}, \S 6, and references therein), are usually
tied to the presentation of Temperley-Lieb algebras. They may,
however, be found in Brauer algebras too (see \cite{W88}, p. 180-181).

According to the isomorphism result for \Kn\ with respect to diagrams
generalizing the friezes of this paper
(see \cite{BDP02} or \cite{DP01}, and references therein; the
first proof of this isomorphism is in
\cite{J83}, \S 4),
we can conclude that the monoid \Dbn\ is isomorphic to a submonoid
of \K2n. Hence we have that \Sn\ is isomorphic to a submonoid of
\K2n. The right-forking term $p^i$ is represented in \K2n\ by
$h^{2i+3}\cirk h^{2i+2}$, and the left-forking term $q^i$
by $h^{2i+1}\cirk h^{2i+2}$.

We can establish also that \Sn\ is isomorphic to submonoids of the
monoids \L2n\ and \J2n\ of \cite{DP01}. These last two monoids
differ from \K2n\ only with respect to the equations $(hc1)$ and
$(hc2)$, which
have to do with circles in diagrams. In \J2n\ we have $c=\mj$,
and instead of $(hc2)$ we have simply $h^i \cirk h^i=h^i$, while in \L2n\
we pay closer attention to circles than in \K2n, as we indicated in Section 6,
and we have something more
involved than $(hc2)$. Circles are, however, irrelevant for the friezes of
this paper, and all of \K2n, \L2n\ and \J2n\ have \Sn\ as a submonoid.

The monoid \Lo\ of \cite{DP01}
has for every $k\in \N^+$ the generators \dk{k} and \gk{k}.
Besides the monoid equations $(1)$ and $(2)$ of Section 2, the following
equations are assumed in \Lo\ for $l\leq k$:

\[
\begin{array}{rl}
\makebox[8em][r]{$({\mbox{\it cup}})$\hspace{1cm}} & \makebox[15em][l]
{$\dk{k}\cirk \dk{l} = \dk{l}\cirk \dk{k\pl 2}$,}
\\[0.1cm]
({\mbox{\it cap}}) \hspace{1cm} & \gk{l}\cirk \gk{k} = \gk{k\pl 2}\cirk \gk{l},
\end{array}
\]
\vspace{-.2cm}
\[
\begin{array}{rl}
\makebox[8em][r]{$({\mbox{\it cup-cap}}\;1)$\hspace{1cm}}
& \makebox[15em][l]
{$\dk{l}\cirk \gk{k\pl 2} =\gk{k}\cirk \dk{l}$,}
\\[0.1cm]
({\mbox{\it cup-cap}}\;2)\hspace{1cm}  & \dk{k\pl 2}\cirk \gk{l}
=\gk{l}\cirk \dk{k},
\\[0.1cm]
({\mbox{\it cup-cap}}\;3)\hspace{1cm}  & \dk{k}\cirk \gk{k\!\pm \!1} = \mj.
\end{array}
\]

\noindent In the monoid \Ko\ we have the additional equation
$\dk{k}\cirk \gk{k} = \dk{l}\cirk \gk{l}$, while in \Jo\ we have
$\dk{k}\cirk \gk{k} = \mj$.
With $h^i$ defined as $\gk{i}\cirk \dk{i}$ and
$c$ as $\dk{i}\cirk \gk{i}$, we can
check easily that \Kn\ is a submonoid of \Ko. It is shown in
\cite{DP01} that \Ko\ is isomorphic to the monoid of
$\cal K$-equivalence classes of the friezes of \cite{DP01},
and analogously for $\cal L$ and $\cal J$
(see Section 6 for $\cal K$-equivalence). These friezes are tied
to self-adjunctions, where an endofunctor is adjoint to itself
(see \cite{DP01}).

Let an order-preserving endomorphism $\varphi$ of \N\ be of {\it
type} $(n,m)$ when for every $k\in \N$ we have that $\varphi(n\pl
k)=m\pl k$. Order-preserving endomorphisms of \N\ that have a type
are said to be of {\it finite type}. It is easy to infer from
Section 6 that for a frieze $D$ of finite type the
order-preserving endomorphism $\varphi(D)$ will be of finite type,
and that for an order-preserving endomorphism $\varphi$ of finite
type the frieze $D(\varphi)$ will be of finite type.

The order-preserving endomorphisms of \N\ of finite type make a monoid
with composition and the identity function on \N. This
monoid, which we call \So, can be presented by generators and relations
in the following manner. The generators of \So\ are for $i\in \N$ the
order-preserving surjective functions $\sigma_i$ such that

\[
\sigma_i(n)=
\left\{
\begin{array}{ll}
n & {\mbox{\rm if }} n\leq i
\\[.1cm]
n\mn 1 & {\mbox{\rm if }} n>i
\end{array}
\right.
\]

\noindent and the order-preserving injective functions $\delta_i$ such
that

\[
\delta_i(n)=
\left\{
\begin{array}{ll}
n & {\mbox{\rm if }} n<i
\\[.1cm]
n\pl 1 & {\mbox{\rm if }} n\geq i.
\end{array}
\right.
\]

The standard generators of the simplicial category \Sc\ are obtained by
restricting the domains and codomains of $\sigma_i$ and $\delta_i$ (see
\cite{ML71}, VII.5). Besides the monoid equations,
these generators satisfy the following
equations for $i\leq j$:

\[
\begin{array}{l}
\sigma_j\cirk \sigma_i=\sigma_i\cirk \sigma_{j+1},\\
\delta_i\cirk \delta_j=\delta_{j+1}\cirk \delta_i,\\
\sigma_i\cirk \delta_{j+2}=\delta_{j+1}\cirk \sigma_i,\\
\sigma_{j+1}\cirk \delta_i=\delta_i\cirk \sigma_j,\\
\sigma_i\cirk \delta_i=\sigma_i\cirk \delta_{i+1}=\mj.
\end{array}
\]

It can be shown that \So\ is isomorphic to a
submonoid of \Ko. The generator
$\sigma_i$ is represented in \Ko\ by $\dk{2i\pl 2}$ and
$\delta_i$ by $\gk{2i\pl 1}$. We can then easily verify the
equations above in \Ko. For the first equation we use
$({\mbox{\it cup}})$, for the second $({\mbox{\it cap}})$,
and for the remaining three $({\mbox{\it cup-cap}}\;1)$,
$({\mbox{\it cup-cap}}\;2)$ and $({\mbox{\it cup-cap}}\;3)$
respectively. We establish in the same manner that \So\ is isomorphic to
submonoids of \Lo\ and \Jo.

The monoids \Sn\ are of course submonoids of \So. We define $p^i$ as
$\delta_{i+1}\cirk \sigma_i$, which is equal to $\sigma_i\cirk \delta_{i+2}$,
and $q^i$ as $\sigma_{i+1}\cirk \delta_i$, which is equal to
$\delta_i\cirk \sigma_i$.

\vspace{.3cm}

\noindent {\footnotesize {\it Acknowledgements.} I would like to
thank Zoran Petri\' c for reading this paper, and an anonymous
referee for informing me about A\u\i zen\v stat's paper, about
which I didn't know when five years ago I wrote the first version
(available at: http://arXiv.org/math.GT/0301302). I would like to
thank also Vitor Fernandes, who was extremely kind to send me a
copy of A\u\i zen\v stat's paper. The writing of the present paper
was financed by the Ministry of Science of Serbia (Grant 144013).}

\end{document}